\documentclass[12pt,english,twoside]{article}
\usepackage[dvips]{graphicx}
\usepackage{authblk}
\usepackage{subeqnarray}
\usepackage{dsfont}
\usepackage{caption}
\usepackage{subcaption}
\usepackage{lipsum}
\usepackage{array,multirow,makecell}
\setcellgapes{1pt}
\RequirePackage{amsmath,amssymb,amsfonts,textcomp,yhmath}
\RequirePackage{times} 
\RequirePackage{cases} 
\usepackage[]{hyperref}
\hypersetup{backref,implicit=true,bookmarks=true, colorlinks=true,
linkcolor=blue,filecolor=blue,
pagecolor=blue,urlcolor=blue,pdfpagemode=UseOutlines}

\newtheorem{theorem}{Theorem}[section]
\newtheorem{lemma}[theorem]{Lemma}
\newtheorem{cor}[theorem]{Corollary}
\newtheorem{proposition}[theorem]{Proposition}

\newtheorem{definition}[theorem]{Definition}

\newtheorem{remark}{Remark}

\newenvironment{preuve}[1][\underline{Proof}]{\textbf{#1.} }{\ \rule{0.5em}{0.5em}}

\begin{document}
\markboth{\small{A. Nangue}}{\tiny{Mathematical analysis of  an
extended cellular model of the Hepatitis C Virus infection with
non-cytolytic process}}

\title{Mathematical analysis of  an extended
cellular model of the Hepatitis C Virus infection with non-cytolytic
process}
\date{}
\author[,1]{Alexis Nangue \thanks{Email : Alexis Nangue : \texttt{alexnanga02@yahoo.fr}}}
\author[2]{Cyprien Fokoue }
\author[3]{Raoue Poumeni \thanks{Now at : elsewhere}}

\affil[1]{University of Maroua, Higher Teacher's Training College,
Department of Mathematics, P.O.Box 55 Maroua, Cameroon}
\affil[2]{University of Maroua, Faculty of Science, Department of
Mathematics and Computer Science, Cameroon}
\affil[3]{University of
Maroua, Higher Teacher's Training College}

\maketitle

\begin{abstract}
The aims of this work is to analyse of the global stability of the
extended model of hepatitis C virus(HCV) infection with cellular
proliferation, spontaneous cure and hepatocyte homeostasis. We first
give general information about hepatitis C. Secondly, We prove the
existence of local, maximal and global solutions of the model and
establish some properties of this solution as positivity and
asymptotic behaviour. Thirdly we show, by the construction of an
appropriate Lyapunov function, that the uninfected equilibrium and
the unique infected equilibrium of the model of HCV are globally
asymptotically stable respectively when the threshold number
$\mathcal{R}_{0}<1-\frac{q}{d_{I}+q}$ and when $\mathcal{R}_{0}>1$.
Finally, some numerical simulations are carried out using Maple
software confirm these theoretical results.
\end{abstract}

{\bf keywords} : HCV model; global solutions; non-cytolytic process;
invariant set; Lyapunov functions; basic reproduction number;
equilibrium points.
\\\\
{\bf AMS Classification Subject 2010} : 92B99, 34D23, 92D25.

\section{Introduction}
Hepatitis C infection is a viral disease caused by Hepatitis C Virus
(HCV) and being transmitted mainly by blood contact between an
infected person and a healthy person. This virus that attacks the
hepatocytes is one of the main causes of chronic diseases of the
liver such as hepatocellular carcinoma, liver cancer and cirrhosis
of the liver \cite{chon, stab}. According to the WHO \cite{oms}
global report published in April 2017 on hepatitis, 200 to 300
million people worldwide are infected with HCV, and between 60$\%$
et 85$\%$ of these people develop chronic liver disease \cite{oms,
anal}. Although there is considerable progress in the research for
the fight against this infection whose virus was discovered in 1989
\cite{chon, stab} and which presents today six genotypes, ranked
from 1 to 6 according to \cite{stab}. There is no vaccine for
prevention yet \cite{relu, oms}. Concerning the treatment of HCV
infection, since 2014, the new direct-acting antivirals combined
with Interferon-$\alpha$ and Ribavirin have been able to cure about
90$\%$ of cases of chronic infection, but leaves the chronic
diseases whose infection has caused. HCV infection is therefore a
major public health problem.
\\\indent
To understand the dynamics of HCV viral load and its infectious
process, mathematical models have become  an important and almost
unavoidable tool\cite{yves, murr}. A model is a system of
mathematical equations accounting for all known experimental data of
the studied biological phenomenon. It makes it possible to better
understand the phenomenon under consideration and to act on the
system optimally. Until 2009, most research work on the modeling of
viral dynamics of HCV only took into account the level of
circulating virus in a human population, the case in vivo was almost
ignored as it provides a better understanding of the pathogenesis of
the virus as Harel Dahari et al \cite{mode} and Chong et
al.\cite{chon}.
\\\indent
Our goal is therefore to analyze the stability of an extended model
of HCV infection in a patient with cell proliferation and
spontaneous healing presented in \cite{Dari05, relu} to reveal
significant information on pathogenesis and dynamics of this virus.
The work is organized as follows : In section 1, first focuses on
presentation of the epidemiological model and give some properties
of its solutions, then we calculate the basic reproduction ratio $
\mathcal{R}_{0}$, which is an indispensable element in the study and
analysis of the models. $ \mathcal{R}_{0}$ is considered in the
virus dynamics as a metric. We theoretically analyze the local
stability and the global stability of the model by the linearization
and Lyapounov's functions respectively in section 2 and in section 3
we perform numerical simulations using biologically plausible
parameter values in Table 1 to confirm the results obtained
theoretically and we complete the work by an appendix.
\section{A dynamical model of HCV infection and some
properties of the solutions of the differential system associated to
the model} Mathematical modeling applies to the study of dynamics of
infectious diseases seems to be one of the most  interesting tool
for designing control or eradication strategies of a disease like
hepatitis C. It allows to test on computer different prevention
scenarios and so help the decision-making in public health
\cite{thes}. In the following sections, we will firstly present
epidemiology and brief history of HCV dynamics, then describe the
model itself, show the existence of solutions and establish some
properties of these solutions and finally calculate the basic
reproduction rate $\mathcal{R}_{0}$ which is an important tool in
epidemiology.
\subsection{Epidemiology and history of HCV}
\subsubsection{HCV infection}
HCV is a virus that attacks the cells of the liver and causes
inflammation of the latter. This virus is present in the blood of an
infected person and is, according to the WHO, mandatory declaration.
It can live for about 5 to 7 weeks in the open air. In the long
term, there may be very serious consequences such as cirrhosis and
in some cases, liver cancer . This virus can remain for decades in
the body without any apparent symptoms. According to \cite{thes,
njou} six main genotypes HCV were identified whereas several
subtypes play an important role in the severity of the disease and
its response to treatment.
\subsubsection{Epidemiology of HCV Infection}
HCV infection is a major public health problem. Worldwide, the
number of chronic HCV carriers is estimated at about
 170 million,
or $ 3 \% $ of the world's population,
and the incidence is between
3 and 4 million infections per year according to the
 statistics
published by the World Health Organization in 2014. The actual
incidence is uncertain   because the distinction between acute and
chronic forms is difficult to make. HCV was only identified in 1989
with the advent of modern molecular cloning techniques.
\\\indent In
Cameroon, about 10,000 people die each year from hepatitis (A, B, C,
D, ...). The country is one of the 17 most affected worldwide, with
a prevalence rate of about $ 13 \% $,  that is 2.5 million people
infected.

\subsubsection{Natural history of HCV infection}
 According to WHO's global report on hepatitis C
 published in April 2017 the first phase of HCV infection
 is said Acute: It can cause
jaundice, but remains asymptomatic in the majority of cases
(70 to
80 $ \% $), hence the risk of going unnoticed.
It is estimated that
20 to 30 $ \% $ people   infected will spontaneously
clear the virus
within the first six months after initial contact.
 If the virus
persists, hepatitis progresses to chronicity.
The liver reacts to
HCV aggression by an inflammatory reaction,
 one of the components of
which is fibrogenesis. Hepatic fibrosis is the main complication of
chronic hepatitis C. Hepatitis C is likely to evolve  at the chronic
phase in about 25 $ \% $ of cases to cirrhosis within a period of 5
to 20 years. In case of cirrhosis, the incidence of hepatocellular
carcinoma is high: on the order of 1 to 4 $ \% $ per year
\cite{thes}.(See figure 1).
\begin{figure}[!h]
\centering
  \includegraphics[width= 14cm, height=5cm]{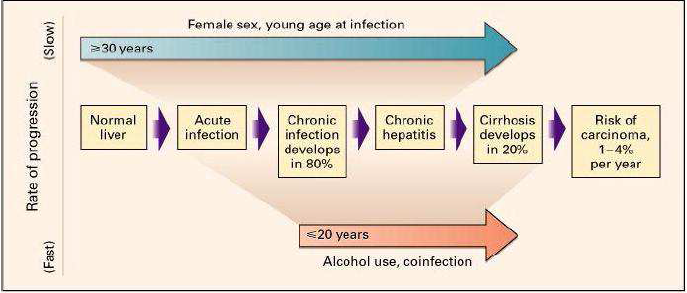}\\
  \caption{Natural history of infection with hepatitis C virus.}
  \label{fig1}
\end{figure}
\\\indent
In order to achieve our various goals, we first describe the model
and its parameters.

\subsection{Description of the HCV model with compartmental
 diagram} 
\begin{figure}[h!]
  \includegraphics[width= 13cm,height=5cm]{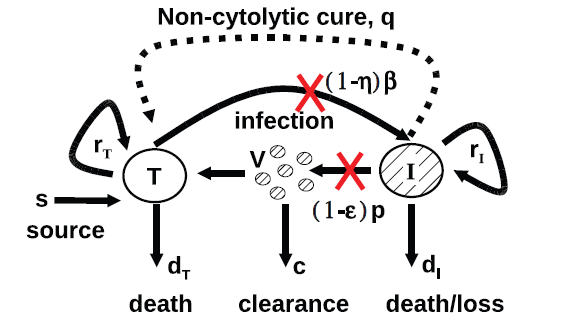}\\
  \caption{Schematic representation of HCV infection models.
   T and I represent target and
infected cells, respectively, and V represents free virus. The
parameters shown in the figure are defined in the text. The original
model of Neumann et al.\cite{Neumann0} assumed that there is no
proliferation of target and infected cells (i.e., $r_{T}=r_{I}=0$)
and no spontaneous cure (i.e., $q = 0$). The extended model of
Dahari, Ribeiro, and Perelson\cite{perel}, which was used for
predicting complex HCV kinetics under therapy, includes target
 and
infected cell proliferation without cure ($r_{T}> 0$, $r_{I} > 0$
and $q = 0$). A model including both proliferation and the
spontaneous cure of infected cells (dashed line; $q > 0$) was used
to explain the kinetics of HCV in primary infection in
chimpanzees\cite{DSRibeiro} }\label{figuure2}
\end{figure}

There are two mathematical models of HCV dynamics :
 the original model or model of Newmann \cite{Neumann0} and its
extended models like that of Dahari \cite{chon, fashion}. Each model
can be represented by a compartmental scheme. A compartmental scheme
is a scheme for estimating the variation in the number of
individuals in each compartment over time. Figure \ref{figuure2} is
the schematic representation of the extended model, which we will
study, of HCV with cellular proliferation and spontaneous healing
designed by T. C. Reluga et al. \cite{relu}. This model expands the
viral dynamics of the original model of infection and the
disappearance of HCV by incorporating the proliferation and death
density dependence. In addition to cell proliferation, the number of
uninfected hepatocytes may increase through immigration or
differentiation of hepatocyte precursors that develop into
hepatocytes at a constitutive rate of $s$ or by spontaneous infected
hepatocyte healing  by a non-cytolytic process at the rate $q$.

The model proposed by Dahari and coworkers \cite{mode, perel}
expands on the standard HCV viral-dynamic model \cite{Neumann0}
 of
infection and clearance by incorporating density-dependent
proliferation and death. Uninfected hepatocytes or noninfected
hepatocytes, T, are infected at a rate $\beta$ per free virus per
hepatocyte. Infected cells, I, produce free virus at rate $p$ per
cell but also die with rate $d_{I}$ . Free virus is cleared at rate
$c$ by immune and other degradation processes. Besides infection
processes, hepatocyte numbers are influenced by homeostatic
processes. Uninfected hepatocytes die at rate $d_{T}$. Both infected
and uninfected hepatocytes proliferate logistically with maximum
rates $r_{I}$ and $r_{T}$ , respectively, as long as the total
number of hepatocytes is less than $T_{max}$. Besides proliferation,
uninfected hepatocytes may increase in number through immigration or
differentiation of hepatocyte precursors that develop into
hepatocytes at constitutive rate $s$, or by spontaneous cure of
infected hepatocytes through a noncytolytic process at rate $q$.
Treatment with antiviral drugs reduces the infection rate by a
fraction $\eta$ and the viral production rate by a fraction
$\varepsilon$. It should be noted that $\eta$ and $\varepsilon$ are
parameters which values are non-negative and less than one.
\\\indent
The interpretations and biologically plausible value of other
parameters are listed in following Table 1 , and a further
comprehensive survey on the description of the model is given in
\cite{Dari05, Dahalay, relu}.
\\\indent
Thus, the variation of healthy hepatocytes T is expressed by the
following expression :
\begin{equation}\label{s1}
 \frac{dT}{dt} = s + r_{T}T \left( 1 -
\frac{T + I}{T_{max}} \right) - d_{T}T - (1- \eta) \beta VT + qI.
\end{equation}
The variation of infected hepatocytes I is expressed by the
following expression:
\begin{equation}\label{s2}
\frac{dI}{dt} = r_{I}I \left( 1 - \frac{T + I}{T_{max}} \right) -
d_{I}I - (1- \eta) \beta VT - qI.
\end{equation}
And the variation of the viral load V is expressed by the following
expression :
\begin{equation}\label{s3}
\frac{dV}{dt} = ( 1 - \varepsilon) pI - cV.
\end{equation}
It follows that the dynamics of T, I and V is governed by the
following differential system :

\begin{equation}\label{s}
\left\lbrace
\begin{array}{lcr}
\dfrac{dT}{dt} = s + r_{T}T \left( 1 - \dfrac{T + I}{T_{max}}
\right) - d_{T}T - (1- \eta) \beta VT + qI
\\\\
\dfrac{dI}{dt} = r_{I}I \left( 1 - d\frac{T + I}{T_{max}} \right) -
d_{I}I - (1- \eta) \beta VT - qI
\\\\
\dfrac{dV}{dt} = ( 1 - \varepsilon) pI - cV
\end{array}
\right. \\
\end{equation}
System (\ref{s}) is under the following initial conditions:
\begin{equation}\label{s0}
T_{0}=T(t_{0}),\quad I_{0}=I(t_{0})\quad \mbox{and} \quad
V_{0}=V(t_{0})\quad \mbox{where} \quad t_{0}\in[0,+\infty[.
\end{equation}

Given the meanings of $\eta$ and $\beta$, the term
 $(1- \eta)\beta VT $ represents the mass action principle; $ \beta VT $ is the
rate of infection of healthy T cells by interaction with
virus V. \\

For biological significance of the parameters, three assumptions are
employed. (1) Due to the burden of supporting virus replication,
infected cells may proliferate more slowly than uninfected cells,
i.e. $r_{I} \leq r_{T}$ . (2) To have a physiologically realistic
model, in an uninfected liver when $T_{max}$ is reached, liver size
should no longer increase, i.e. $s \leq d_{T}T_{max} $. (3) Infected
cells have a higher turnover rate than uninfected cells, i.e. $d_{I}
\geq d_{T}$. The interpretations and biologically plausible values
of other parameters are listed in Table 1, and a further
comprehensive survey on the description of (\ref{s}) is given in
\cite{relu}. Besides HCV infection, the similar model of (\ref{s})
is also used to describe the dynamics of HBV or HIV infection, in
which the full logistic terms mean the proliferation of
uninfected/infected hepatocytes \cite{SMCiu, QLi, DSRibeiro}, or the
mitotic transmission of uninfected/infected CD4+T.
 \\
The range of variation of each parameter is recorded in table 1
\cite{relu}.
 \begin{center}
 \underline{Table 1}
 \end{center}
 \small{Estimated parameter ranges for
 hepatitis C when modeled with system \cite{relu}. The $r_{T}$,
  $T_{max}$, and $d_{T}$ parameters are not independently
identifiable, so common practice is to fix $d_{T}$ prior
 to fitting}
\label{tb}
$$\begin{tabular}{|cccc|}
  \hline
 Symbol & Minimum & Maximum & Units \\\hline
  $\beta$  & $10^{-8}$ & $10^{-6}$ & $virus^{-1}.ml.day^{-1}$ \\
  $T_{max}$ & $4 \times 10^{6}$ & $1.3.10^{7}$ & $cells.ml^{-1}$ \\
  $p$& $0.1$ & $44$ & $virus.cell^{-1}.day^{-1}$\\
  $s$ & $1$ & $1.8 \times10^{5}$ & $cells.ml^{-1}.day^{-1}$ \\
  $q$ & $0$ & $1$ & $day^{-1}$\\
  $c$ & $0.8$ & $22$ & $day^{-1}$\\
  $d_{T}$ & $10^{-3}$ & $1.4 \times 10^{-2}$ & $day^{-1}$\\
  $d_{I}$ & $10^{-3}$ & $0.5$ & $day^{-1}$\\
  $r_{T}$ & $2\times 10^{-3}$ & $3.4$ & $day^{-1}$\\
  $r_{I}$ & $Unknown$ & $Unknown$ & $day^{-1}$ \\
  \hline
\end{tabular}$$
This table tells us in which interval varies each parameter of the
model. For the study of stability and for simulations these ranges
of values will have to be respected for a good decision-making.
\subsection{Theorems of existence and some properties of solution
to the cauchy problem (\ref{s}), (\ref{s0})}
\subsubsection{Existence of local solutions}
\begin{theorem}\label{t1}
Let $ T_{0}, I_{0}, V_{0} \in \mathbb{R}$. There exists $t_{1} > 0$
and functions $T, I, V : [t_{0};  t_{1} [ \longrightarrow
\mathbb{R}$ continuously differentiable such that $(T, I, V)$ is a
solution of system (\ref{s}) satisfying (\ref{s0}).
\end{theorem}
\begin{preuve}
We will use the local Cauchy-Lipschitz theorem
to proof this. Since the system of equations (\ref{s})
is autonomous, it is enough to show that the function \\
$$\begin{tabular}{c c c l}
$f$ :&$\mathbb{R}^{3}$&$\longrightarrow$ &$\mathbb{R}^{3}$\\
&$(T, I, V)$&$\longmapsto$&$\Big{(}f_{1}(T, I, V); f_{2}(T, I, V);
f_{3}(T, I, V)\Big{)}$
\end{tabular}$$
is locally Lipschitzian with:
\begin{equation}\label{f1}
  f_{1}(T, I, V) = s+
r_{T}T\left(1-\frac{T+I}{T_{max}}\right) - d_{T}T - (1- \eta) \beta
V T+ qI.
\end{equation}
\begin{equation}\label{f2}
f_{2}(T, I, V) = r_{I}I\left(1-\frac{T+I}{T_{max}}\right) + (1-
\eta) \beta V T - d_{I}I- qI.
\end{equation}
\begin{equation}\label{f3}
 f_{3}(T, I, V) = (1- \varepsilon)pI -cV.
\end{equation}
  According to L. Perko \cite{lawr}, it is also enough to prove that
   $f$ is a class $\mathcal{C}^{1}$ function.\\ The jacobian matrix of
    $f$ at
$(T, I, V)$ is: \\
$$ Df(T, I, V) =
\begin{pmatrix}
\frac{\partial f_{1}}{\partial T}(T, I, V) & \frac{\partial
f_{1}}{\partial I}(T, I, V) & \frac{\partial f_{1}}{\partial V}(T,
I,
V)\\\\
\frac{\partial f_{2}}{\partial T}(T, I, V) &\frac{\partial
f_{2}}{\partial I}(T, I, V) & \frac{\partial f_{2}}{\partial V}(T,
I,
V)\\\\
\frac{\partial f_{3}}{\partial T}(T, I, V) &\frac{\partial
f_{3}}{\partial I}(T, I, V) & \frac{\partial f_{3}}{\partial
V}(T,I,V)
\end{pmatrix},$$ \\
i.e.
\\ $$ Df(T, I, V)   =
\begin{pmatrix} r_{T}\left(1-\frac{2T+I}{T_{max}}\right) - d_{T}+ &
-\frac{r_{T}T}{T_{max}} + q & - (1- \eta) \beta
T\\
- (1- \eta) \beta V&&\\\\
 -\frac{r_{I}I}{T_{max}}+ (1- \eta) \beta V &
r_{I}\left(1-\frac{T+2I}{T_{max}}\right) - d_{I} - q & (1- \eta)
\beta
T\\\\
0 & (1- \varepsilon)p & - c
\end{pmatrix}.$$
Each component of this matrix being continuous, they are locally
bounded for all $ (T, I, V) \in \mathbb{R}^{3}$. Therefore $f$
possesses continuous and bounded partial derivatives on
 any compact
of $\mathbb{R}$. Thus  $f$ is locally Lipschitzian with respect to
$(T, I, V)$. By the Cauchy-Lipschitz theorem, there is a local
solution defined on  $[t_{0}; t_{1}[$ This completes the proof of
this theorem~\ref{t1}.
\end{preuve}
\begin{remark}
The function $f$ in the proof of theorem~\ref{t1} is a class
 $C^{1}$ function
so the system (\ref{s}) has a unique maximal solution.
\end{remark}

\subsubsection{positivity of the  system (\ref{s})}
\begin{theorem}\label{t2}
Let $(T, I, V)$ be a solution of the system (\ref{s}) over an
interval $[t_{0}, t_{1}[$
 such that $ T(t_{0})=T_{0}, I(t_{0})=I_{0}$ et $V(t_{0}) = V_{0}$.\\
  If  $T_{0}, I_{0}, V_{0}$ are positive, then $T(t)$, $I(t)$ and $V(t)$
  are also positive for all $ t \in
[t_{0}, t_{1}[$.
\end{theorem}
\begin{preuve} We are going to prove by contradiction. so
suppose there is  $ t \in [t_{0}, t_{1}[$ such that $T(t)=0$
or $I(t)=0$ or $V(t)=0$.\\
Let $x=(x_{1}, x_{2}, x_{3})=(T, I, V)$\\
Let also $t_{*}$ be the smallest of all $t$ in the interval $[t_{0},
t_{1}[$ such that $x_{i}(t) > 0,\\ \forall t \in [t_{0}, t_{*}[,
\forall i \in \left\lbrace 1, 2, 3 \right\rbrace $ and
 $x_{i}(t_{*})
=0$ for a certain $i$.\\ Then each of the equations of the
 system
(\ref{s}) can be written $\dot{x}_{i} = - h_{i}(x) +
 g_{i}(x)$ where $g_{i}$ is a non negative function and $h_{i}$
 any function.\\
Thus,
\begin{eqnarray*}
\quad\dot{T} & = & \frac{dT}{dt}, \\
  & = & - T \left( -r_{T} \left( 1 - \frac{T+I}{T_{max}} \right) + d_{T} + (1- \eta ) \beta V \right) + s +
  qI,
  \\
  & = & - T h_{1}(T, I,  V) + g_{1}(T, I, V).
\end{eqnarray*}
with
$$h_{1}(T, I,  V)=\left( -r_{T} \left( 1 - \frac{T+I}{T_{max}}
\right) + d_{T} + (1- \eta ) \beta V \right)$$
and
$$g_{1}(T, I,
V)=s + qI ;$$
similarly;
 \begin{eqnarray*}
\quad\dot{I} & = & \frac{dI}{dt}, \\
  & = & - I \left( -r_{I} \left( 1 - \frac{T+I}{T_{max}} \right) + d_{I} + q  \right) + (1- \eta ) \beta V
  T,
  \\
  & = & - I h_{2}(T, I,  V) + g_{2}(T, I, V).\\
\end{eqnarray*}
with $$h_{2}(T, I,  V)=\left( -r_{I} \left( 1 - \frac{T+I}{T_{max}}
\right) + d_{I} + q  \right)$$ and
 $$g_{2}(T, I, V)=(1- \eta ) \beta
V T$$ and
\begin{eqnarray*}
V & = & \frac{dV}{dt}, \\
  & = & - Vc + ( 1- \varepsilon)pI ,\\
  & = & -V h_{3}(T, I,  V) + g_{3}(T, I, V).
\end{eqnarray*}
with
 $$h_{3}(T, I,  V)=c$$
 and
 $$g_{3}(T, I, V)=( 1-
\varepsilon)pI .$$
Without loss the generality, suppose that $x_{1}(t_{*}) =0$.\\
As hypothesized, $g_{1}(T, I, V)$ is positive on $[t_{0}, t_{*}]$,
it follows that
$$ \dot{T} \geq - T h_{1}(T, I,  V);$$  from where
 $$\frac{d}{dt}(\log T) \geq - h_{1}(T, I,  V).$$
  Yet $(T, I, V) $
is a solution of (\ref{s}). $T$, $I$, $V$ are class class
 $C^{1}$ functions.
They are continuous on $[t_{0}, t_{*}]$
 and therefore are bounded on $[t_{0}, t_{*}]$.
  Thus $h_{1}(T, I,  V)$ is bounded on  $[t_{0}, t_{*}]$.\\
 There exists a constant  $ k > 0$ such that $$\frac{d}{dt}(\log T) \geq - h_{1}(T, I,  V) > - k.$$
 By integrating this previous expression  on  $[t_{0}, t_{*}]$,
 we get  $$\log T(t_{*}) - \log T(t_{0})
 \geq - k( t_{*} - t_{0}).$$ Where $$ T(t_{*}) \geq T(t_{0})
   e^{-k (t_{*}-t_{0})} > 0.$$ This is a
   contradiction because $ T(t_{*})=0$.\\\indent
Similarly, assuming that $x_{2}(t_{*}) =0$, as hypothesized,
$g_{2}(T, I, V)$ is positive on  $[t_{0}, t_{*}]$,
it follows that
 $$ \dot{I} \geq - I h_{2}(T, I,  V).$$ Thus,
$$\frac{d}{dt}(\log I) \geq - h_{1}(T, I,  V).$$
Yet $(T, I, V) $ is a solution of  system (\ref{s}); $T$, $I$, $V$
are class $ C^{1}$ functions. They are therefore continuous on
  $[t_{0}, t_{*}]$
 and consequently are bounded on $[t_{0}, t_{*}]$.\\
 Therefore  $h_{2}(T, I,  V)$ is bounded on $[t_{0}, t_{*}]$.\\
 Thus, there exists a constant $ \lambda > 0$ such that
  $$\frac{d}{dt}(\log I) \geq - h_{2}(T, I,  V) > - \lambda.$$
  Integrating, the last expression on  $[t_{0}, t_{*}]$, yields
  $$\log I(t_{*}) - \log I(t_{0}) \geq - \lambda( t_{*} - t_{0}).$$
  Therefore $$ I(t_{*}) \geq I(t_{0})
   e^{-\lambda (t_{*}-t_{0})} > 0.$$ This is a contradiction because
    $ I(t_{*})=0$.\\\indent
As far as that goes, assuming that $x_{3}(t_{*}) =0$.\\
By hypothesis, $g_{3}(T, I, V)$ is positive on $[t_{0}, t_{*}]$, we
obtain $V \geq - Vc.$ i.e
 $$\frac{d}{dt}(\log I) \geq -c.$$ Hence integrating on
  $[t_{0}, t_{*}]$ we obtain $$ V(t_{*}) \geq V(t_{0})e^{-k
  (t_{*}-t_{0})} > 0.$$ This is a contraction.\\
Conclusion: $ T$, $I$, $V$ are positive on  $[t_{0}, t_{1}[$.
\end{preuve}\\
It will now be shown, with the help of the continuation criterion
the existence of global solutions of problem (\ref{s}), (\ref{s0}).
\subsubsection{Existence of global solutions}
\begin{theorem}
The solutions  of the Cauchy problem (\ref{s}), (\ref{s0}),
 with
positive initial data,
 exist globally in time
in the future that is  on $[t_{0}, +\infty[$.
\end{theorem}
\begin{preuve}
To prove this it is enough to show that all variables
 are bounded on
an arbitrary finite interval $[t_{0}; t)$. Using the positivity, by
the theorem~\ref{t2}, of the solutions it is enough to show that all
variables
 are bounded above.\\\indent
 Taking the sum
of  equations (\ref{s1}), (\ref{s2}) shows that :
$$ \frac{d}{dt}(T+I) \leq \lambda $$
 and hence that $T(t)+I(t)\leq  T_{0}+I_{0}+\lambda(t-t_{0})$.
 Thus $T$ and $I$ are
bounded on any finite interval. The third equation, i.e.
$$\frac{dV}{dt} = ( 1 - \varepsilon) pI - cV,$$
 then shows that $V(t)$
cannot grow faster than linearly and is also bounded
 on any finite interval.
This completes the proof of this theorem.
\end{preuve}
\subsubsection{Asymptotic behaviour}
\begin{theorem}\label{t3}
For any positive solution $(T, I, V)$ of system (\ref{s}),
(\ref{s0}) we have :
\begin{equation*}
    T(t) \leq \tilde{T}_{0} \;\;,I(t) \leq \tilde{T}_{0}\;\;  and \;\;
    V(t)\leq \lambda_{0}
\end{equation*}
where
\begin{equation*}
    \tilde{T}_{0}=\frac{T_{max}}{2r_{I}} \left( \sqrt{(r_{T} -
d_{T})^{2} + \frac{4sr_{I}}{T_{max}}} + r_{T} - d_{T} \right),
    \lambda_{0}=\lambda_{0} = \max \Big{\{} V_{0},
        \frac{( 1 - \varepsilon)}{c}p\tilde{T_{0}}\Big{\}}.
\end{equation*}
\end{theorem}
\begin{preuve}
Summing  equations (\ref{s1}) and (\ref{s2}), we get:
\begin{eqnarray*}
\frac{d}{dt}(T+I)  & = & s + \left( 1 - \frac{T+I}{T_{max}} \right) \left( r_{T}T - r_{I}I \right) - d_{T}T - d_{I}I, \\
      & = & s + r_{T}T - r_{I}I - d_{T}T - d_{I}I - \frac{T+I}{T_{max}}, \\
      & = & s + \left( r_{T} - d_{T} \right) T + \left( r_{I} - d_{I} \right)I - \frac{T+I}{T_{max}}(r_{T}T + r_{I}I), \\
      & \leqslant & s + \left( r_{T} - d_{T} \right)( T + I ) - \frac{T+I}{T_{max}} (r_{T}T + r_{I}I) \quad\ \hbox{since}
      \quad \ r_{T} - d_{T} \geq r_{I} - d_{I},\\
\mbox{thus}\; \frac{d }{dt}(T+I)  & \leqslant & s + \left( r_{T} -
d_{T} \right)( T + I ) - \frac{r_{I}}{T_{max}} (T + I)^{2}\quad \
\hbox{since} \quad r_{I} \leqslant
      r_{T}.
\end{eqnarray*}
Let $N_{1} = T + I$,\quad $a = s > 0$,\quad $b = \left( r_{T} -
d_{T} \right) > 0$,\quad $d = - \frac{r_{I}}{T_{max}} < 0$ and let
us solve the following equation
\begin{equation}\label{n1}
\frac{dN_{1}}{dt} = a + bN_{1} + dN_{1}^{2}
\end{equation}
Coupled to  equation (\ref{n1}) the initial condition :
\begin{equation}\label{n2}
N_{1}(t_{0}) = N_{1}^{0}.
\end{equation}
The resolution of the problem (\ref{n1}), (\ref{n2}) gives for all
$t\in[t_{0},+\infty[$,
 $$ N_{1}(t) = - \frac{1}{2d}\biggl[ \tanh \biggl(
\frac{1}{2}\sqrt{-4ad + b^{2}} - \frac{1}{2}t_{0} \sqrt{-4ad +
b^{2}} - \arctan \left( \frac{2N_{1}^{0} + b}{\sqrt{-4ad + b^{2}}}
\right) \biggr) \sqrt{-4ad + b^{2}} \biggr] - \frac{b}{2d}.$$ As for
all $ x \in \mathbb{R}, -1 \leqslant \tanh x \leqslant 1$,
 it follows that :
\begin{equation*}
  N_{1}(t) \leqslant -
\frac{1}{2d} \left( \sqrt{-4ad + b^{2}} + b \right)
\end{equation*}
i.e.
\begin{equation*}
 N_{1}(t)
\leqslant \frac{T_{max}}{2r_{I}} \left( \sqrt{(r_{T} - d_{T})^{2} +
\frac{4sr_{I}}{T_{max}}} + r_{T} - d_{T} \right).
\end{equation*}
Let $$ \tilde{T}_{0} = \frac{T_{max}}{2r_{I}} \left( \sqrt{(r_{T} -
d_{T})^{2} + \frac{4sr_{I}}{T_{max}}} + r_{T} - d_{T} \right), $$ we
obtain :
$$ N_{1}(t) \leqslant \tilde{T}_{0}.$$
Therefore $$ T +I \leqslant \tilde{T}_{0}.$$ Since T and I are
positive $ I \leqslant T + I$ and $ T \leqslant T + I$, so it
follows that $T(t)\leqslant \tilde{T}_{0} $ and $I(t)\leqslant
\tilde{T}_{0} $.\\
From (\ref{s3}), we have:
\begin{eqnarray*}
 \frac{dV}{dt} & \leqslant & ( 1 - \varepsilon)p(T + I) - cV, \\
               &  \leqslant & ( 1 - \varepsilon)p\tilde{T_{0}} - cV \quad\ \hbox{since}\quad \ T +I \leqslant
               \tilde{T_{0}}.
\end{eqnarray*}
According to Gronwall inequality,
\begin{eqnarray*}
  V(t) & \leqslant & V(t_{0}) e^{-c (t-t_{0})} + \int_{t_{0}}^{t} ( 1 - \varepsilon)p\tilde{T_{0}} e^{\int_{u}^{t} - c ds} du, \\
       & \leqslant & V_{0} e^{-c (t-t_{0})} + ( 1 - \varepsilon)p\tilde{T_{0}} \int_{t_{0}}^{t} e^{-c(t-u)} du,
\end{eqnarray*}
\begin{eqnarray*}
  V(t)  & \leqslant & V_{0} e^{-c (t-t_{0})} + ( 1 - \varepsilon)p\tilde{T_{0}} \frac{e^{-c (t-t)}-e^{-c (t-t_{0})}}{c},
\end{eqnarray*}
\begin{eqnarray*}
     & \leqslant & V_{0} e^{-c (t-t_{0})} + ( 1 - \varepsilon)p\tilde{T_{0}} \frac{1-e^{-c (t-t_{0})}}{c}, \\
       & \leqslant & \lambda_{0} \Big{(} e^{-c (t-t_{0})} + 1-e^{-c (t-t_{0})}\Big{)} \\
 V(t) & \leqslant & \lambda_{0}.
\end{eqnarray*}
with $$ \lambda_{0} = \max \Big{\{} V_{0},
        \frac{( 1 - \varepsilon)}{c}p\tilde{T_{0}}\Big{\}}.$$
This completes the proof of theorem~\ref{t3}.
\end{preuve}
\begin{remark}
 It follows that all solutions of the system (\ref{s}) are asymptotically
uniformly bounded in compact subset  $\Omega$ defined by
\begin{equation*}
 \Omega = \left\lbrace (T, I, V) \in \mathbb{R} ;
0< T + I \leqslant\tilde{T_{0}} ;\quad 0< V \leqslant \lambda_{0}
\right\rbrace .
\end{equation*}
\end{remark}
\begin{remark}
Theorem~\ref{t3} shows that all solutions of model (\ref{s}) in
$\mathds{R}^{3}_{+}$ are ultimately bounded and according to
theorem~\ref{t2}, that solutions with positive initial value
conditions are positive, which indicates that model (\ref{s}) is
well-posed and biologically valid.
\end{remark}
\subsection{Basic reproduction ratio $\mathcal{R}_{0}$}
One of the most important concerns about any infectious disease is
its ability to invade a population. Many epidemiological models
have
a disease free equilibrium (DFE) at which the population
remains in
the absence of disease. These models usually have a threshold
parameter, known as the basic reproduction number,
$\mathcal{R}_{0}$, such that if $\mathcal{R}_{0} < 1$,
 then the DFE
is locally asymptotically stable, and the disease
cannot invade the
population, but if $\mathcal{R}_{0} > 1$,
 then the DFE is unstable
and invasion is always possible. In other words, we have the
following definition :
\begin{definition}\cite{ODiek}
The basic reproduction ratio or the basic reproduction number or
basic reproductive ratio  $\mathcal{R}_{0}$ is defined as the
expected number of secondary cases produced, in a completely
susceptible population, by a typical infected individual during its
entire period of infectiousness.
\end{definition}
Determine $ \mathcal{R}_{0} $ in function of  the parameters of the
model allow us to guess the conditions under which the disease
invade the population.
\subsubsection{Determination of
the uninfected equilibrium or virus-free equilibrium or noninfected
equilibrium}
\begin{proposition}\label{t0}
The uninfected equilibrium point $E^{0}$ of the system (\ref{s}) is
given by
$$E^{0}=(T^{0},0,0)$$
 where :
 $$T^{0} =
\frac{T_{max}}{2r_{T}} \left(r_{T} - d_{T} + \sqrt{(r_{T} -
d_{T})^{2} + \frac{4r_{T}s}{T_{max}}} \right).$$
\end{proposition}
\begin{preuve}
When there is no viral infection, the uninfected hepatocytes dynamic
is determined by :
\begin{equation}
\frac{dT}{dt} = s + r_{T}T \left( 1- \frac{T}{T_{max}} \right)  -
d_{T}T.
\end{equation}
The quantity $T^{0}$ of the free-virus equilibrium point $E^{0}$
is solution of the equation  $\frac{dT}{dt}=0$.\\
Hence, let us solve the following equation :
$$-\frac{r_{T}}{T_{max}} T^{2} + (r_{T} - d_{T}) T + s =0.$$
Its discriminant is given by :
 $$ \Delta = (r_{T} - d_{T})^{2} + 4 \frac{r_{T}}{T_{max}}s ,$$
 which yields :
 $$T= \frac{d_{T}-r_{T} - \sqrt{(r_{T} - d_{T})^{2} + 4
\frac{r_{T}}{T_{max}}s}}{2\frac{r_{T}}{T_{max}}}.$$ Thus, in the
absence of viral infection, the amount of susceptible cells or
uninfected hepatocytes attend to a positive constant level $T^{0}$,
which is :
$$T^{0} =
\frac{T_{max}}{2r_{T}} \left(r_{T} - d_{T} + \sqrt{(r_{T} -
d_{T})^{2} + 4 \frac{r_{T}}{T_{max}}s} \right). $$ This completes
the proof.
\end{preuve}
\subsubsection{Computation of the basic reproduction number $\mathcal{R}_{0}$}
We are going to use  Van Den Driessche and Watmough method
\cite{ODiek, PVan} for calculating the basic reproduction ratio
 $\mathcal{R}_{0}$ of the model (\ref{s}).\\\indent
 Let us first present briefly the method.
 \\\indent
Considering population whose population are grouped into $n$
homogeneous compartments $ X= (X_{1}, X_{2}, ..., X_{n})$ where
$X_{i} \geq 0$  is the number of individuals in compartment  $i$.
For clarity we sort the compartments $X_{i}$, $i= 1, ..., n$
 so that
the first $m$ $( m \leqslant n )$ compartments correspond to
infected individuals. The distinction between infected and
uninfected compartments must be determined from the epidemiological
interpretation of the model and cannot be deduced from the structure
of the equations alone.\\
We define $X_{s}$ to be the set of all disease free states. That is
: $$ X_{s} = \left\lbrace X \geq 0 / X_{1}= X_{2}= ....= X_{m}=0
 \right\rbrace.$$
Let $F_{i}(X)$ be the rate of appearance of new infections in
compartment $i$ that is the infected individuals coming from  other
compartments and enter
into $i$.\\
$V_{i}^{+}$ be the rate of transfer of individuals into compartment
$i$ by all other means (displacement, healing, aging).\\
$V_{i}^{-}$  be the rate of transfer of individuals out of
compartment $i$ (mortality, change of statut).\\
It is assumed that each function is continuously
differentiable at
least twice in each variable.
\\\indent
Figure 3 below shows the variations of the number of
  individuals in compartment $i$ in a population.
 \begin{figure}[h!]
 \begin{center}
  \includegraphics[width= 5cm,height=4cm]{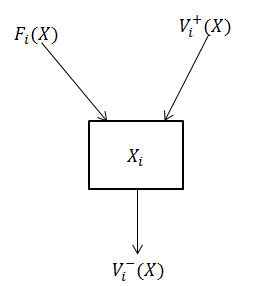}\\
  \caption{\textbf{variation of the number of
  individuals in compartment $i$ in a population}}
  \end{center}
\end{figure}
The variation of the number of individuals in compartment $i$ is
given by :
$$ \frac{dX_{i}}{dt}= F_{i}(X) - V_{i}(X)$$
where
$$V_{i}(X) = V_{i}^{-}(X) - V_{i}^{+}(X).$$

Due to the nature of the epidemiological model, we have the
following properties :
 \begin{enumerate}
\item [$P_{1})$] If $X_{i} \geq 0$ then  $F_{i}(X) \geq 0$,
 $V_{i}^{-}(X) \geq 0$ and $V_{i}^{+}(X) \geq 0$.\\
 Since
each function represents a directed transfer of individuals.
\item [$P_{2})$] If $  X_{i} = 0$ then $ V_{i}^{-}(X)=0$.\\
Indeed, If a compartment is empty, then there can be no
 transfer of
individuals out of the compartment by death, infection,
 nor any
other means :  it is the essential property of a compartmental
model.
\item [$P_{3})$] Pour $ i > m$, $F_{i}(X) = 0$. \\
In fact, the compartments with an index greater than $m$ are
"uninfected". By definition, it can not appear in these compartments
infected individual.
\item [$P_{4})$] Si $X \in X_{s}$ alors $F_{i}(X) = 0$ et pour
 $i \leqslant m$,
$V_{i}^{+}(X) = 0$.\\
Indeed, to ensure that the disease free subspace is invariant, we
assume that if the population is free of disease then the population
will remain free of disease. That is, there is no (density
independent) immigration of infectives. This is Lavoisier's
principle. There is no spontaneous generation.
\end{enumerate}
Let $ F = (F_{i})_{i=1,...,n}$ and
$V = V^{+} - V^{-} = (V_{i}^{+} - V_{i}^{-})_{i = 1,.....,n}$.\\
Let also $x^{0}$ the uninfected equilibrium point of the
corresponding model.\\
Let $DF(x^{0})$ and $DV(x^{0})$  denote the jacobian matrices of F
and V respectively at the point $x^{0}$.\\
It follows that $ DF(x^{0})$ is a positive matrix and
 $DV(x^{0})$ a Metzler matrix  (matrix whose the
extra diagonal terms are greater  or equal than zero ). Thus we have
the following equivalent definition of $\mathcal{R}_{0}$ :
\begin{definition}\cite{ODiek, PVan}
$$\mathcal{R}_{0} = \rho (- DF.DV^{-1})$$
 where $\rho$ represents the spectral radius. i.e
$$\mathcal{R}_{0} = \max \vert \lambda \vert$$ with $$ \lambda \in
sp ( - DF.DV^{-1}),$$ where $sp(A)$ is the spectrum of  $A$, i.e.
the set of eigenvalues associated to a matrix $A.$
\end{definition}
\subsubsection{Expression of the basic reproduction number
 $\mathcal{R}_{0}$ associated to the system (\ref{s})}
\begin{proposition}
 The expression of the basic reproduction number
 $\mathcal{R}_{0}$ associated to the system (\ref{s})
 is given by :
\begin{equation}\label{r0}
\mathcal{R}_{0} = \frac{r_{I}}{d_{I}+q } \left( 1-
\frac{T^{0}}{T_{max}} \right) + \frac{(1-\theta)\beta T^{0}p}
{c(d_{I}+q )}.
\end{equation}
where  $$1-\theta=(1- \varepsilon )(1- \eta ).$$
\end{proposition}
\begin{preuve}
Concerning the model (\ref{s}) that we study here, the system of
infected states is the following :
\begin{center}
 $\left\lbrace
\begin{array}{lcr}
\dfrac{dI}{dt} = r_{I}I \left( 1- \dfrac{T+I}{T_{max}} \right) + (1-
\eta ) \beta VT - d_{I}I -
qI\\\\
\dfrac{dV}{dt} = (1- \varepsilon )pI -cV
\end{array}
\right. $
\end{center}
The expression of the quantities $F$, $V$, $DF(E^{0})$, $DV(E^{0})$
and $DV^{-1}(E^{0})$ are given by: \\
 $ F=
\begin{pmatrix}
 r_{I}I \left( 1- \frac{T+I}{T_{max}}\right)
  + (1- \eta ) \beta VT
 \\\\
0
\end{pmatrix}$,\quad
$V=
\begin{pmatrix}
(- d_{I} - q )I \\\\
(1- \varepsilon )pI -cV
\end{pmatrix}$,
\\
 $DF(E^{0})=
\begin{pmatrix}
r_{I} \left( 1- \frac{T^{0}}{T_{max}} \right) \quad & (1- \eta )
\beta VT^{0}
\\\\
0 \quad & 0
\end{pmatrix}$,
$DV(E^{0})=
\begin{pmatrix}
- d_{I} - q \quad & 0 \\\\
(1- \varepsilon )p\quad & -c
\end{pmatrix}$,
\\
$DV^{-1} =
\begin{pmatrix}
\frac{-1}{d_{I}+ q}\quad & 0 \\\\
\frac{-(1- \varepsilon )p}{c( d_{I} + q )}\quad & - \frac{1}{c}
\end{pmatrix}.$
\\ \\
From those quantities, we obtain:\\

$DF.DV^{-1}=
\begin{pmatrix}
-r_{I} \left( 1- \frac{T^{\circ}}{T_{max}} \right)\quad & -(1- \eta)
\beta
VT^{0}\\\\
0 \quad & 0
\end{pmatrix}$
$\begin{pmatrix}
\frac{-1}{d_{I}+ q}\quad & 0 \\\\
\frac{-(1- \varepsilon )p}{c( d_{I} + q )}\quad & - \frac{1}{c}
\end{pmatrix}$ \\\\

\hspace{2cm} = $\begin{pmatrix} \frac{r_{I}}{d_{I} + q } \left( 1-
\frac{T^{0}}{T_{max}} \right) + \frac{(1- \varepsilon )(1- \eta
)\beta T^{0}p}{c(d_{I}+ q)}\quad & \frac{(1- \eta ) \beta T^{0}}{c}
\\\\
0\quad & 0
\end{pmatrix}$.

Let   $1-\theta=(1- \varepsilon )(1- \eta )$. It follows that :
\begin{eqnarray*}
| - DF.DV^{-1} - \lambda I_{2}|=0 & \Leftrightarrow & -\lambda
\left( \frac{r_{I}}d_{I}+q   \left( 1- \frac{T^{0}}{T_{max}} \right)
 + \frac{(1-\theta)\beta T^{0}p}{c(d_{I}+q )} - \lambda \right)
 =0\\
 & \Leftrightarrow & \lambda = 0 \;\; \mbox{ou} \;\;
 \lambda = \frac{r_{I}}d_{I}+q  \left( 1- \frac{T^{0}}{T_{max}} \right)
 + \frac{(1-\theta)\beta T^{0}p}{c(d_{I}+q )}.
\end{eqnarray*}
Therefore :
\begin{equation*}
\mathcal{R}_{0} = \frac{r_{I}}{d_{I}+q  } \left( 1-
\frac{T^{0}}{T_{max}} \right) + \frac{(1-\theta)\beta
T^{0}p}{c(d_{I}+q )},\\
\end{equation*}
which completes the proof of the proposition.
\end{preuve}
\begin{remark}
$\theta\in]0, 1[$ denotes the overall effectiveness rate of the drug
.
\end{remark}
\begin{remark}
Henceforth, we will let $\delta = d_{I}+q$ and $1-\theta=(1-
\varepsilon )(1- \eta )$.
\end{remark}
At the end of this section, we note that HCV is a major
health problem in the world and particularly in Cameroon where it
affects almost 13\% of population.\\\indent
 For the model (\ref{s}) which is the
subject of our work, we have shown the existence of the global
solution
 and
establish some properties like positivity. The calculation of
 $\mathcal{R}_{0}$ has been done.
\\\\\indent In the
next section, we will determine the  infected equilibrium point
 and
 establish
the conditions on $\mathcal{R}_{0}$ for which  stability of the
model occurs.

\section{Stability analysis of the model}
In this section, we study the local stability and global stability
of the equilibrium points and we present some numerical
 simulations
of the theoretical results obtained. Specifically, we prove by
Lyapunov's theory that the uninfected equilibrium point $ E^{0}$
 is
globally asymptotically stable if
$\mathcal{R}_{0}<1-\frac{q}{\delta}$ and the infected equilibrium
point $ E^{\ast} $ is globally asymptotically stable when it exists.
\\\indent Before that we establish a number of essential preliminary
results for the next steps
\subsection{Invariant set of the model}
\begin{theorem}\label{invset}
Let $(t_{0}, S_{0}=(T_{0}, I_{0}, V_{0})) \in \mathds{R}\times
\mathds{R}^{3}_{+} $ and $([t_{0}, T[,  S=(T, I, V))$ be a maximal
solution of the Cauchy problem (\ref{s}), (\ref{s0})  ($T\in ]t_{0},
+\infty[$). If $T(t_{0})+I(t_{0})\leq \tilde{T}_{0}$ and
$V(t_{0})\leq \lambda_{0}$ then the set :
  $$ \Omega = \left\lbrace (T, I, V) \in \mathbb{R} ;
0 < T + I \leqslant\tilde{T}_{0} ;\quad 0 < V \leqslant \lambda_{0}
\right\rbrace ,$$
  where :
\begin{equation*}
    \tilde{T}_{0}=\frac{T_{max}}{2r_{I}} \left( \sqrt{(r_{T} -
d_{T})^{2} + \frac{4sr_{I}}{T_{max}}} + r_{T} - d_{T}
\right),\;\;\mbox{and}\;\;
    \lambda_{0} \max \Big{\{} V_{0},
        \frac{( 1 - \varepsilon)}{c}p\tilde{T}_{0}\Big{\}},
\end{equation*}
is a positively invariant set by system (\ref{s}).
\end{theorem}
\begin{preuve}
Let $t_{1} \in [t_{0}, T[$. We shall show that :
\begin{enumerate}
    \item[(i)]If $T(t_{1})+I(t_{1})\leq \tilde{T}_{0} $ then for all
    $t_{1}\leq t < T$, $T(t)+I(t)\leq \tilde{T}_{0}$.
    \item[(ii)]If $V(t_{1})\leq \lambda_{0}$ then for all $t_{1}\leq t < T$,  $V(t)\leq
    \lambda_{0}$.
\end{enumerate}
\begin{enumerate}
    \item[i)] Let us show i) by contradiction.\\
Let us suppose that there exists  $\varepsilon_{1}>0$ such that
$t_{1}<t_{1}+\varepsilon_{1}<+\infty$ we have
$$(T+I)(t_{1}+\varepsilon_{1})>T_{0}.$$
Let $t_{1}^{\ast}=inf\{t\geq
t_{1};(T+I)(t_{1}+\varepsilon_{1})>T_{0}\}.$\\
 If
$(T+I)(t_{1}^{\ast})=T_{0},$ then since
$$(T+I)(t)=T_{0} +
\frac{d}{dt}\left(T(t_{1}^{\ast})+I(t_{1}^{\ast})
\right)(t-t_{1}^{\ast})+o(t-t_{1}^{\ast})$$  when $t\rightarrow
t_{1}^{\ast}$. In addition, according to equations (\ref{s1}) and
(\ref{s2}) of system (\ref{s}), we have
$$\frac{d}{dt}\left((T+I)(t_{1}^{\ast})\right)\leq
s+(r_{T}-d_{T})T_{0}-\frac{r_{I}}{T_{max}}T_{0}^{2}.$$ Recall that
$$s+(r_{T}-d_{T})T_{0}-\frac{r_{I}}{T_{max}}T_{0}^{2}=0.$$
It follows that :
$$\frac{d}{dt}\left((T+I)(t_{1}^{\ast})\right)\leq0.$$
Hence, there exists $\tilde{\varepsilon}>0$ such that for all
$t\in[t_{1}^{\ast},t_{1}^{\ast}+\varepsilon[$, $(T+I)(t)\leq T_{0},$
which is a contradiction. Therefore for all  $t\in [t_{0},+\infty[$,
$(T+I)(t)\leq T_{0}$.

\item[ii)]Let us show ii) by contradiction.\\
Let us suppose that there exists $\varepsilon_{1}>0$ such that
$t_{1}<t_{1}+\varepsilon_{1}<+\infty$  and
$$T(t_{1}+\varepsilon_{1})>\lambda_{0}.$$
 Let
$t_{2}^{\ast}=inf\{t\geq
t_{1};V(t)>\lambda_{0}\}$. \\
Since $V(t_{2}^{\ast})= \lambda_{0}$, since:
$$V(t)=\lambda_{0}+\frac{dV}{dt}(t_{2}^{\ast})(t-t_{2}^{\ast})
+o(t-t_{2}^{\ast})$$ with  when $t\rightarrow
t_{2}^{\ast}$.\\
Equation (\ref{s3}) yields :
\begin{eqnarray*}
 \frac{d}{dt}(V(t_{2}^{\ast}))& \leq
  &(1-\varepsilon)p(I+T)(t_{2}^{\ast})-cV(t_{2}^{\ast})\\
                                      & \leq &
                                      (1-\varepsilon)pT_{0}
                                      -c\lambda_{0};
                                      \quad
\end{eqnarray*}
yet
$$\lambda_{0}=\max\left\{V_{0},
\frac{(1-\varepsilon)pT_{0}}{c}\right\};$$
 consequently
 \begin{eqnarray*}
 \frac{d}{dt}(V(t_{2}^{\ast})) & \leq & (1-\varepsilon)pT_{0}-c\frac{(1-\varepsilon)pT_{0}}{c}\\
                                      &\leq &0.
\end{eqnarray*}
Thus, there exists $\tilde{\varepsilon}>0$ such that for all
$t_{2}^{\ast}\leq t\leq t_{2}^{\ast}+\varepsilon$,  $V(t)\leq
\lambda_{0}$ which is a contradiction. Therefore for all
$t\in[t_{0},+\infty[$,\quad $V(t)\leq \lambda_{0}.$ Which completes
the proof of Theorem~\ref{invset}.
\end{enumerate}
\end{preuve}

\subsection{ Existence of the infected equilibrium point}
When it exists, the infected equilibrium point is given by:
$E^{\ast}=(T^{\ast},I^{\ast},V^{\ast})$ where $T^{\ast}$, $I^{\ast}$
and  $V^{\ast}$ are positive constants that we are going to
determine.
\begin{lemma}\label{lem211}
$T^{\ast}$ exists if and only if
$$s+q\frac{T_{max}}{r_{I}}(r_{I}-\delta)>0.$$
\end{lemma}
\begin{preuve}
Let us consider the following system of algebraic equations : \\
\begin{numcases}\strut
s + r_{T}T\left(1 - \frac{T + I}{T_{max}}\right) - d_{T}T - (1-
\eta) \beta VT + qI=0; \label{se1}\\
 r_{I}I\left(1 - \frac{T +
I}{T_{max}}\right) - d_{I}I - (1- \eta) \beta VT - qI=0;
\label{se2}\\
 ( 1 - \varepsilon) pI - cV=0. \label{se3}
\end{numcases}
  (\ref{se3}) yields :
  \begin{equation}\label{se4}
V=\frac{(1-\varepsilon)pI}{c}.
\end{equation}
Reporting (\ref{se4}) in (\ref{se2}),  we have :
\\\vspace{0.5cm}
$$r_{I}I\left(1-\frac{T+I}{T_{max}}\right)+\frac{(1-\eta)(1-\varepsilon)\beta
pIT}{c}-\delta I=0.$$
Hence
\begin{equation*}
r_{I}\left(1-\frac{T+I}{T_{max}}\right)+\frac{(1-\theta)\beta
pT}{c}-\delta =0\;\; \mbox{since} \;\; I\neq 0 \;\; \mbox{(there is
an infection)}
\end{equation*}
i.e
\begin{equation*}
 \frac{r_{I}I}{T_{max}}=\frac{(1-\theta)\beta pT}{c}-\delta
r_{I}\left(1-\frac{T}{T_{max}}\right).
\end{equation*}
It follows that :
\begin{equation}\label{se5}
I=\left(\frac{(1-\theta)\beta
pT_{max}}{cr_{I}}-1\right)T+\frac{T_{max}}{r_{I}}(r_{I}-\delta).
\end{equation}
Let
\begin{equation}\label{se6}
 g_{1}(T)=qI=q \left(\frac{(1-\theta)\beta
pT_{max}}{cr_{I}}-1 \right)T+q\frac{T_{max}}{r_{I}}(r_{I}-\delta).
\end{equation}
Reporting  (\ref{se4}) and (\ref{se5}) in (\ref{se1}) leads to  :
\begin{eqnarray*}
   && s+r_{T}T\left(1-\frac{T}{T_{max}}\right)-\frac{r_{T}T}{T_{max}}\left(\left(\frac{(1-\theta)\beta
pT_{max}}{cr_{I}}-1\right)T+\frac{T_{max}}{r_{I}}(r_{I}-\delta)\right)-d_{T}T
  \\
   && -(1-\eta)\beta
T\frac{(1-\varepsilon)p}{c}\left(\left(\frac{(1-\theta)\beta
pT_{max}}{cr_{I}}-1\right)T+\frac{T_{max}}{r_{I}}(r_{I}-\delta)\right)\\
   &&+q\left(\left(\frac{(1-\theta)\beta
pT_{max}}{cr_{I}}-1\right)T+\frac{T_{max}}{r_{I}}(r_{I}-\delta)\right)
= 0
\end{eqnarray*}
i.e.
\begin{eqnarray*}
   && s+T\left(r_{T}-d_{T}-\frac{r_{T}}{r_{I}}(r_{I}-\delta)-\frac{(1-\theta)\beta
pT_{max}}{cr_{I}}(r_{I}-\delta)\right)-\frac{(1-\theta)\beta
p}{c}\Bigl(\frac{r_{T}}{r_{I}}+\frac{(1-\theta)\beta
pT_{max}}{cr_{I}}-1\Bigr)T^{2}\\
   &&+q\left(\frac{(1-\theta)\beta
pT_{max}}{cr_{I}}-1\right)T+q\frac{T_{max}}{r_{I}}(r_{I}-\delta)=0.
\end{eqnarray*}
 Thus ,
\begin{eqnarray*}
   && s+q\frac{T_{max}}{r_{I}}(r_{I}-\delta)+T\biggl(r_{T}-d_{T}-\frac{r_{T}}{r_{I}}(r_{I}-\delta)-\frac{(1-\theta)\beta
pT_{max}}{cr_{I}}(r_{I}-\delta)+q\left(\frac{(1-\theta)\beta
pT_{max}}{cr_{I}}-1\right)\biggr)  \\
   &&-\frac{(1-\theta)\beta
p}{c}\left(\frac{r_{T}}{r_{I}}+\frac{(1-\theta)\beta
pT_{max}}{cr_{I}}-1\right)T^{2}=0.
\end{eqnarray*}

It follows that,
\begin{eqnarray*}
   && s+q\frac{T_{max}}{r_{I}}(r_{I}-\delta)+T\Bigl(r_{T}-\delta-\frac{r_{T}}{r_{I}}(r_{I}-\delta)-\frac{(1-\theta)\beta
pT_{max}}{cr_{I}}(r_{I}-\delta) +\frac{(1-\theta)\beta
pT_{max}}{cr_{I}}q\Bigr) \\
   &&-\frac{(1-\theta)\beta
p}{c}\left(\frac{r_{T}}{r_{I}}+\frac{(1-\theta)\beta
pT_{max}}{cr_{I}}-1\right)T^{2}=0.
\end{eqnarray*}
Let
\begin{eqnarray*}
 h_{2}(T)&=&s+q\frac{T_{max}}{r_{I}}(r_{I}-\delta)+T\Bigl(r_{T}-\delta-\frac{r_{T}}{r_{I}}(r_{I}-\delta)-\frac{(1-\theta)\beta
pT_{max}}{cr_{I}}(r_{I}-\delta)+\frac{(1-\theta)\beta
pT_{max}}{cr_{I}}q\Bigr)  \\
   &&-\frac{(1-\theta)\beta
p}{c}\left(\frac{r_{T}}{r_{I}}+\frac{(1-\theta)\beta
pT_{max}}{cr_{I}}-1\right)T^{2};
\end{eqnarray*}
 we have :
\begin{equation}\label{h2}
h_{2}(T)=g_{2}(T)+g_{1}(T)
\end{equation}
with :
\begin{equation*}
g_{2}(T)=s+T\left(r_{T}-d_{T}-\frac{r_{T}}{r_{I}}(r_{I}-\delta)-\frac{(1-\theta)\beta
pT_{max}}{cr_{I}}(r_{I}-\delta)\right)-\frac{(1-\theta)\beta
p}{c}\left(\frac{r_{T}}{r_{I}}+\frac{(1-\theta)\beta
pT_{max}}{cr_{I}}-1\right)T^{2}.
\end{equation*}
Since $r_{I}\leq r_{T}$, the polynomial (\ref{h2}) has a unique
positive root $T^{\ast}$ if and only if :
$$s+q\frac{T_{max}}{r_{I}}(r_{I}-\delta)>0.$$
This completes the proof of Lemma~\ref{lem211}.
\end{preuve}\\
Suppose that $\delta \geq r_{I}$, we have the following results :
\begin{lemma}\label{lem222}
If $\frac{(1-\theta)\beta pT_{max}}{cr_{I}}>1$, then :
\begin{enumerate}
    \item[i)] $T^{\ast}\leq
\tilde{T_{1}}$ if and only if  $g_{1}(T^{\ast})\leq 0$
    \item[ii)] $T^{\ast}> \tilde{T_{1}}$ if and only if
     $g_{1}(T^{\ast})> 0$.
\end{enumerate}
\end{lemma}
\begin{remark}
$\tilde{T}_{1}$ of Lemma~\ref{lem222} is the solution of equation
$g_{1}(T)= 0,$ i.e
$$\tilde{T}_{1}=\frac{c(\delta-r_{I})T_{max}}{(1-\theta)\beta
pT_{max}-cr_{I}}.$$
\end{remark}
\begin{proposition}\label{p3}
Suppose that  $\delta\geq r_{I}$.
\begin{itemize}
    \item
    If
$\frac{(1-\theta)\beta pT_{max}}{cr_{I}}\leq 1$, then
$g_{1}(T^{\ast})\leq 0$. Hence, system (\ref{s})  admits no infected
equilibrium point.
    \item
    If $\frac{(1-\theta)\beta
pT_{max}}{cr_{I}}>1$ then :
     \begin{enumerate}
        \item[i)]system (\ref{s})  admits no infected
equilibrium point when $T^{\ast}\leq \tilde{T_{1}}$.
        \item[ii)] system
(\ref{s}) admits a unique infected equilibrium point $E^{\ast}$ when
$T^{\ast}> \tilde{T_{1}}$.
     \end{enumerate}
\end{itemize}
\end{proposition}
\begin{lemma}\label{la}
$\tilde{T_{1}}=\frac{(r_{I}-\delta)T^{0}}{r_{I}-\delta
\mathcal{R}_{0}}$.
\end{lemma}
\begin{preuve}
Since $\tilde{T}_{1}$ is the root of equation (\ref{se6}), we have :
\begin{eqnarray*}
 \tilde{T_{1}}& = &\frac{c(r_{I}-\delta)T_{max}}{cr_{I}-(1-\theta)}\beta pT_{max}\\
                                      & = & (r_{I}-\delta)\left(\frac{1}{\frac{r_{I}}{T_{max}}-\frac{(1-\theta)\beta p}{c}}\right)\\
                                      & = &(r_{I}-\delta)\left(\frac{T^{0}}{\frac{r_{I}T^{0}}{T_{max}}-\frac{(1-\theta)\beta pT^{0}}{c}}\right)\\
                                      & = &\frac{(r_{I}-\delta)T^{0}}{r_{I}-\delta
                                     \mathcal{R}_{0}}.
\end{eqnarray*}
\end{preuve}
\begin{lemma}\label{lb}
$\frac{(1-\theta)\beta pT_{max}}{cr_{I}}\leq1$ if and only if
$\mathcal{R}_{0}\leq\frac{r_{I}}{\delta}$.
\end{lemma}
 \begin{preuve} We have :
\begin{eqnarray*}
 \frac{(1-\theta)\beta pT_{max}}{cr_{I}}\leq 1 & \Leftrightarrow &\frac{(1-\theta)\beta p}{c}\leq\frac{r_{I}}{T_{max}}\\
                                      & \Leftrightarrow &\frac{(1-\theta)\beta pT^{0}}{c\delta}\leq\frac{r_{I}T^{0}}{\delta T_{max}} \\
                                      & \Leftrightarrow &\frac{(1-\theta)\beta pT^{0}}{c\delta}-\frac{r_{I}T^{0}}{\delta T_{max}}+\frac{r_{I}}{\delta}\leq\frac{r_{I}}{\delta}\\
                                      &\Leftrightarrow & \mathcal{R}_{0}\leq \frac{r_{I}}{\delta}.
\end{eqnarray*}
\end{preuve}

Since $\frac{r_{I}}{\delta}\leq 1$, the equivalence of the
Lemma~\ref{lb} allows us to write  $\mathcal{R}_{0}\leq1$.\\
The following proposition establishes the link between
$\mathcal{R}_{0}$ and the existence of the equilibrium point
$E^{\ast}$ when  $\frac{(1-\theta)\beta pT_{max}}{cr_{I}}>1$.
\begin{proposition}\label{prop3.7}
Suppose that  $T^{\ast}$ exists, $r_{I}\leq\delta$ and
$\frac{(1-\theta)\beta pT_{max}}{cr_{I}} > 1$, then :
\begin{enumerate}
    \item[i)]  $T^{\ast}>
\tilde{T_{1}}$ if $\mathcal{R}_{0}> 1$.
    \item[ii)] $T^{\ast} \leq
\tilde{T_{1}}$ if $\frac{r_{I}}{\delta}<\mathcal{R}_{0}\leq1$.
\end{enumerate}
\end{proposition}
\begin{preuve}
Recall that $T^{\ast}> \tilde{T_{1}}$ if and only if
$h_{2}(\tilde{T_{1}})>0$ and  $T^{\ast}\leq \tilde{T_{1}}$ if and
only if  $h_{2}(\tilde{T_{1}})\leq0.$\\
Using the expression of  $h_{2}(T)$ given in (\ref{h2}) we have :
$$h_{2}(\tilde{T_{1}})=g_{2}(\tilde{T_{1}}).$$
Thus
\begin{eqnarray*}
h_{2}(\tilde{T_{1}}) &=& s+\frac{(\delta-r_{I})T^{0}}{\delta
\mathcal{R}_{0}-r_{I}}\left(r_{T}-d_{T}-\frac{r_{T}}{r_{I}}(r_{I}-\delta)-\frac{(1-\theta)\beta
pT_{max}}{cr_{I}}(r_{I}-\delta)\right)  \\
   && -\frac{(1-\theta)\beta
p}{c}\left(\frac{r_{T}}{r_{I}}+\frac{(1-\theta)\beta
pT_{max}}{cr_{I}}-1\right)\left(\frac{(\delta-r_{I})T^{0}}{\delta
\mathcal{R}_{0}-r_{I}}\right)^{2}
\end{eqnarray*}
 Yet
\begin{eqnarray*}
 s &=& \left(d_{T}-r_{T}\left(1-\frac{T^{0}}{T_{max}}\right)\right)T^{0},\\
 \frac{(1-\theta)\beta pT_{max}}{cr_{I}}  &=& \frac{T_{max}}{T^{0}}\left(\frac{\delta
 \mathcal{R}_{0}-r_{I}}{r_{I}}\right)+1
\end{eqnarray*}
and
\begin{eqnarray*}
 \frac{(1-\theta)\beta p}{c}  &=& \frac{\delta
\mathcal{R}_{0}-r_{I}}{T^{0}}+\frac{r_{I}}{T_{max}}.
\end{eqnarray*}
 Hence,
\begin{eqnarray*}
h_{2}(\tilde{T_{1}})&=&\biggl(d_{T}-r_{T}(1-\frac{T^{0}}{T_{max}})\biggr)T^{0}+\biggl(r_{T}-d_{T}-\frac{r_{T}}{r_{I}}(r_{I}-\delta)+\\\\
&&-\frac{T_{max}}{T^{0}}\Bigl(\frac{\delta
\mathcal{R}_{0}-r_{I}}{r_{I}}\Bigr)(r_{I}-\delta)-(r_{I}-\delta)\biggr)\frac{(\delta-r_{I})T^{0}}{\delta
\mathcal{R}_{0}-r_{I}}+\\\\
&&-\Bigl(\frac{\delta
\mathcal{R}_{0}-r_{I}}{T^{0}}+\frac{r_{I}}{T_{max}}\Bigr)\biggl(\frac{r_{T}}{r_{I}}+\frac{T_{max}}{T^{0}}\Bigl(\frac{\delta
\mathcal{R}_{0}-r_{I}}{r_{I}}\Bigr)
\biggr)\Bigl(\frac{{(}\delta-r_{I}{)}T^{0}}{\delta
\mathcal{R}_{0}-r_{I}}\Bigr)^{2}\\\\
&=&\frac{\delta T^{0}(\mathcal{R}_{0}-1)}{(\delta
\mathcal{R}_{0}-r_{I})^{2}T_{max}}\Biggl(r_{T}T^{0}(\delta-r_{I})+(\delta
\mathcal{R}_{0}-r_{I})T_{max}\biggl(d_{T}-r_{T}\Bigl(1-\frac{T^{0}}{T_{max}}\Bigr)\biggr)
\Biggr).
\end{eqnarray*}
Since $d_{T}-r_{T}\left(1-\frac{T^{0}}{T_{max}}\right)>0$ and
$\delta\geq r_{I}$, we get :
 $$h_{2}(\tilde{T_{1}})>0\quad
 \mbox{if} \quad \mathcal{R}_{0}>1$$
 and
$$h_{2}(\tilde{T_{1}})\leq 0\quad
\mbox{if} \quad
 \frac{r_{I}}{\delta}< \mathcal{R}_{0}\leq
1.$$ This completes the proof of Proposition~\ref{prop3.7}.
\end{preuve}\\
Suppose now  $\delta < r_{I}$, then equation (\ref{h2}) admits a
unique positive solution   $T^{\ast}$ and we have the following
results:
\begin{lemma}\label{lem333}
Suppose that $\frac{(1-\theta)\beta pT_{max}}{cr_{I}}<1$, then :
\begin{itemize}
    \item $T^{\ast}\geq \tilde{T_{2}}$ if and only if
     $g_{1}(T^{\ast})\leq 0$.
    \item $T^{\ast}< \tilde{T_{2}}$ if and only if $g_{1}(T^{\ast})> 0$.
\end{itemize}
\end{lemma}
\begin{remark}
$\tilde{T_{2}}$ of lemma~\ref{lem333}is the solution of equation
$g_{1}(T)= 0$, i.e
$$\tilde{T_{2}}=\frac{c(r_{I}-\delta)T_{max}}{cr_{I}-(1-\theta)\beta
pT_{max}}.$$
\end{remark}
\begin{proposition}\label{p2}
\begin{enumerate}
    \item[i)] If $\frac{(1-\theta)\beta pT_{max}}{cr_{I}}\geq 1$,
     then
$g_{1}(T)>0$ and the system (\ref{s}) admits in this case a unique
infected equilibrium point  $E^{\ast}$.
    \item[ii)] If $\frac{(1-\theta)\beta
pT_{max}}{cr_{I}}<1$, alors :
                \begin{itemize}
                         \item the system (\ref{s}) not admits an
                         infected equilibrium point when $T^{\ast}\geq \tilde{T_{2}}$.
                         \item the system (\ref{s}) admits a
unique infected equilibrium point $E^{\ast}$ when $T^{\ast}<
\tilde{T_{2}}$.
                \end{itemize}
\end{enumerate}
\end{proposition}
We state the following two lemmas whose the proofs are analogous of
those of lemma~\ref{la} and  lemma~\ref{lb} respectively. These
lemmas will help us to complete the conditions of existence of the
infected equilibrium point $E^{\ast}.$
\begin{lemma}
$\tilde{T_{2}}=\frac{(r_{I}-\delta)T^{0}}{r_{I}-\delta
\mathcal{R}_{0}}$
\end{lemma}
\begin{lemma}
$\frac{(1-\theta)\beta pT_{max}}{cr_{I}}\geq1$ if and only if
$\mathcal{R}_{0}\geq\frac{r_{I}}{\delta}$.
 \end{lemma}
\begin{proposition}\label{p1}
Suppose that : $r_{I}>\delta$ and $\frac{(1-\theta)\beta
pT_{max}}{cr_{I}} < 1$.\\
Then $T^{\ast} \geq \tilde{T_{2}}$ if $\mathcal{R}_{0}\leq 1$, and
$T^{\ast} < \tilde{T_{2}}$ if
$1<\mathcal{R}_{0}<\frac{r_{I}}{\delta}$.
\end{proposition}
\begin{preuve}
Recall that $T^{\ast}\geq \tilde{T_{2}}$ if and only if
$h_{2}(\tilde{T_{2}})\geq 0$ and $T^{\ast}< \tilde{T_{2}}$ if and
only if  $h_{2}(\tilde{T_{2}})<0.$ From the expression of $h_{2}(T)$
given by (\ref{h2}) we have :
            $$h_{2}(\tilde{T_{2}})=g_{2}(\tilde{T_{2}}).$$
Thus :
\begin{eqnarray*}
h_{2}(\tilde{T_{2}})   &=& s+\frac{(r_{I}-\delta)T^{0}}{r_{I}-\delta
\mathcal{R}_{0}}\left(r_{T}-d_{T}-\frac{r_{T}}{r_{I}}(\delta-r_{I})-\frac{(1-\theta)\beta
pT_{max}}{cr_{I}}(\delta-r_{I})\right) \\
                        && -\frac{(1-\theta)\beta
p}{c}\left(\frac{r_{T}}{r_{I}}+\frac{(1-\theta)\beta
pT_{max}}{cr_{I}}-1\right)\left(\frac{(r_{I}-\delta)T^{0}}{r_{I}-\delta
\mathcal{R}_{0}}\right)^{2}
\end{eqnarray*}
 yet
$$s=\left[d_{T}-r_{T}\left(1-\frac{T^{0}}{T_{max}}\right)\right]T^{0},$$
 $$\frac{(1-\theta)\beta pT_{max}}{cr_{I}}=\frac{T_{max}}{T^{0}}
 \left(\frac{r_{I}-\delta \mathcal{R}_{0}}{r_{I}}\right)+1$$
$$  \mbox{and} \quad \frac{(1-\theta)\beta p}{c}=\frac{r_{I}-\delta
\mathcal{R}_{0}}{T^{0}}+\frac{r_{I}}{T_{max}}.$$

Hence :
\begin{eqnarray*}
h_{2}(\tilde{T_{2}})&=&\left(d_{T}-r_{T}\left(1-\frac{T^{0}}{T_{max}}\right)\right)T^{0}+\biggl(r_{T}-d_{T}-\frac{r_{T}}{r_{I}}(\delta-r_{I})+\\\\
&&-\frac{T_{max}}{T^{0}}\left(\frac{r_{I}-\delta
\mathcal{R}_{0}}{r_{I}}\right)(\delta-r_{I})-(\delta-r_{I})\biggr)\frac{(r_{I}-\delta)T^{0}}{r_{I}-\delta
\mathcal{R}_{0}}+\\\\
&&-\left(\frac{r_{I}-\delta
\mathcal{R}_{0}}{T^{0}}+\frac{r_{I}}{T_{max}}\right)\left(\frac{r_{T}}{r_{I}}+\frac{T_{max}}{T^{0}}\left(\frac{r_{I}-\delta
\mathcal{R}_{0}}{r_{I}}\right)
\right)\left(\frac{(r_{I}-\delta)T^{0}}{r_{I}-\delta
\mathcal{R}_{0}}\right)^{2}\\\\
&=&\frac{\delta T^{0}(1-\mathcal{R}_{0})}{(r_{I}-\delta
\mathcal{R}_{0})^{2}T_{max}}\left[r_{T}T^{0}(r_{I}-\delta)+(r_{I}-\delta
\mathcal{R}_{0})T_{max}\left(d_{T}-r_{T}\left(1-\frac{T^{0}}
{T_{max}}\right)\right) \right].
\end{eqnarray*}
Since $d_{T}-r_{T}\Big{(}1-\frac{T^{0}}{T_{max}}\Big{)}>0$ and
$\delta < r_{I}$, we get :
\begin{itemize}
    \item  $h_{2}(\tilde{T_{2}})\geq0$
 if $\mathcal{R}_{0}\leq 1$.
    \item $h_{2}(\tilde{T_{2}})< 0$ if
$1<\mathcal{R}_{0}<\frac{r_{I}}{\delta}$.
\end{itemize}
This completes the the proof.
\end{preuve}\\
The conditions of existence of an infected equilibrium point
$E^{\ast}$ have been established, we are going in the following
subsection give its expression.
\subsection{Expression of the equilibrium points}
\subsubsection{Uninfected equilibrium point}
By the  proposition~\ref{t0}, the uninfected equilibrium point is
given by  $E^{0}=(T^{0},0,0)$ where
$$T^{0}=\frac{T_{max}}{2r_{T}}\left(r_{T}-d_{T}+\sqrt{(r_{T}-d_{T})^{2}+\frac{4sr_{T}}{T_{max}}}\right).$$
\subsubsection{Infected equilibrium point}
\begin{lemma}\label{l1}
When it exists, $T^{\ast}$ is defined by :
$$T^{\ast}=\frac{1}{2}\left(-\frac{D}{H}+\sqrt{\left(\frac{D}{H}\right)^{2}+F+\frac{4sT_{max}}{r_{T}}H}\right)$$
where :
$$D=AT_{max}\left(\frac{1}{r_{T}}\Big{(}1+\frac{d_{T}+q}{A}\Big{)}-\frac{\delta}{r_{I}}\Big{(}\frac{1}{r_{T}}+\frac{1}{A}\Big{)}-\frac{q}{r_{T}r_{I}}\right);$$
$$F=\frac{4AqT_{max}^{2}}{H^{2}r_{T}^{2}r_{I}^{2}}\Big{(}A(\delta-r_{I})-d_{I}(r_{I}-r_{T})-r_{I}(q-r_{I}-r_{T})+r_{T}q\Big{)};$$
$$H=\frac{A^{2}}{r_{I}r_{T}}+\frac{A}{r_{I}}-\frac{A}{r_{T}}$$
 and
$$A=\frac{(1-\theta)\beta pT_{max}}{c}.$$
\end{lemma}
\begin{preuve}
 When $T^{\ast}$  exists, it will be a positive solution of the
 equation of second degree :
 \begin{equation}\label{h}
h_{2}(T)=0,
\end{equation}
with :
\begin{eqnarray*}
h_{2}(T)&=&s+q\frac{T_{max}}{r_{I}}(r_{I}-\delta)+T\Bigl(r_{T}-\delta-\frac{r_{T}}{r_{I}}(r_{I}-\delta)-\frac{(1-\theta)\beta
pT_{max}}{cr_{I}}(r_{I}-\delta)+\frac{(1-\theta)\beta
pT_{max}}{cr_{I}}q\Bigr)  \\
   &&-\frac{(1-\theta)\beta
p}{c}\left(\frac{r_{T}}{r_{I}}+\frac{(1-\theta)\beta
pT_{max}}{cr_{I}}-1\right)T^{2}.
\end{eqnarray*}
Let
$$A=\frac{(1-\theta)\beta pT_{max}}{c};$$ then (\ref{h})
becomes :
\begin{equation*}
s+q\frac{T_{max}}{r_{I}}(r_{I}-\delta)+\left(r_{T}-(d_{T}+q)
-\frac{r_{T}}{r_{I}}(r_{I}-\delta)-\frac{A}{r_{I}}
\left(r_{I}-\delta +\frac{Aq}{r_{I}}\right)\right)T-
\frac{r_{T}}{T_{max}} \left(\frac{A^{2}}{r_{I}r_{T}}
+\frac{A}{r_{I}}-\frac{A}{r_{T}}\right)T^{2}=0.
\end{equation*}
Let also :
$$H=\frac{A^{2}}{r_{I}r_{T}}+\frac{A}{r_{I}}-\frac{A}{r_{T}},$$
then (\ref{h}) becomes :
\begin{equation*}
s+q\frac{T_{max}}{r_{I}}(r_{I}-\delta)+\left(r_{T}-(d_{T}+q)-\frac{r_{T}}{r_{I}}(r_{I}-\delta)-\frac{A}{r_{I}}\left(r_{I}-\delta+
\frac{Aq}{r_{I}}\right)\right)T-\frac{r_{T}}{T_{max}}HT^{2}=0.
\end{equation*}
Let :
\begin{equation*}
   a=-\frac{r_{T}H}{T_{max}}      ,
\end{equation*}
\begin{equation*}
  d=s+q\frac{T_{max}}{r_{I}}(r_{I}-\delta)      ,
\end{equation*}
\begin{equation*}
b=r_{T}-(d_{T}+q)-\frac{r_{T}}{r_{I}}(r_{I}-\delta)-
\frac{A}{r_{I}}\Big{(}r_{I}-\delta+\frac{Aq}{r_{I}}\Big{)} ,
\end{equation*}
hence the previous equation yields :
$$aT^{2}+bT+d=0.$$
Its discriminant is:
$$\Delta=b^{2}-4ad$$ i.e.
$$\Delta=\left[r_{T}-(d_{T}+q)-\frac{r_{T}}{r_{I}}\left(r_{I}-\delta-\frac{A}{r_{I}}(r_{I}-\delta)+\frac{Aq}{r_{I}}\right)\right]^{2}
+\left(s+\frac{qT_{max}}{r_{I}}\delta\right).$$
Hence
\begin{eqnarray*}
 T^{\ast} & = &\frac{r_{T}-(d_{T}+q)-\frac{r_{T}}{r_{I}}(r_{I}-\delta-\frac{A}{r_{I}}(r_{I}-\delta)+\frac{Aq}{r_{I}}+\sqrt{\Delta}}{2\frac{Hr_{T}}{T_{max}}}\\
                                      & = & \frac{1}{2}\left(\frac{AT_{max}\left(\frac{r_{T}-(d_{T}+q)}{Ar_{T}}-\frac{r_{T}}{r_{I}r_{T}A}(r_{I}-\delta)
                                      -(\frac{1}{r_{T}}-\frac{\delta}{r_{I}r_{T}})+\frac{q}{r_{I}r_{T}}\right)}{H}+\frac{T_{max}\sqrt{\Delta}}{Hr_{T}}\right) \\
\end{eqnarray*}
\begin{eqnarray*}
                                      & = &\frac{1}{2}\left(\frac{AT_{max}\left(\frac{r_{T}-(d_{T}+q)}{Ar_{I}}-\frac{1}{r_{I}A}(r_{I}-\delta)
                                      -\Big{(}\frac{1}{r_{T}}-\frac{\delta}{r_{I}r_{T}}\Big{)}+\frac{q}{r_{I}r_{T}}\right)}{H}+
                                      \frac{T_{max}\sqrt{\Delta}}{Hr_{T}}\right)\\
                                      & = &\frac{1}{2}\left(\frac{AT_{max}\left(\frac{1}{r_{T}}(1+\frac{d_{T}+q}{A}-\frac{\delta}{r_{I}}(\frac{1}{r_{T}}+
                                      \frac{1}{A})-\frac{q}{r_{I}r_{T}}\right)}{H}+\frac{T_{max}\sqrt{\Delta}}{Hr_{T}}\right)\\
                                      & = &
                                      \frac{1}{2}\left(-\frac{D}{H}+\frac{T_{max}\sqrt{\Delta}}{Hr_{T}}\right).
\end{eqnarray*}
Where
$$D=AT_{max}\left(\frac{1}{r_{T}}\Big{(}
1+\frac{d_{T}+q}{A}\Big{)}-\frac{\delta}{r_{I}}
\Big{(}\frac{1}{r_{T}}+\frac{1}{A}\Big{)}-
\frac{q}{r_{T}r_{I}}\right).$$
 We have :
\begin{eqnarray*}
 \frac{T_{max}\sqrt{\Delta}}{Hr_{T}} & =
 &\sqrt{\frac{T_{max}^{2}\Delta}{H^{2}r_{T}^{2}}}\\\\
                                      & = & \sqrt{\left(\frac{D}{H}\right)^{2}-\frac{4Hr_{T}T_{max}}{H^{2}r_{T}^{2}}\left(s+
                                      \frac{qT_{max}}{r_{I}}(r_{I}-\delta)\right)}
                                      \\\\
                                      & = &\sqrt{\left(\frac{D}{H}\right)^{2}+\frac{4sT_{max}}{Hr_{T}}+
                                      \frac{4AqT_{max}^{2}}{H^{2}r_{T}^{2}r_{I}^{2}}\left(\frac{Hr_{T}r_{I}^{2}}{A}-\frac{Hr_{T}r_{I}}{A}\delta\right)}
                                      \\\\
                                      & = &\sqrt{\left(\frac{D}{H}\right)^{2}+\frac{4sT_{max}}{Hr_{T}}+
                                      \frac{4AqT_{max}^{2}}{H^{2}r_{T}^{2}r_{I}^{2}}\left(A(\delta-r_{I})-d_{I}(r_{I}-r_{T})\right)-r_{I}(q-r_{I}-r_{T})+r_{T}q}
                                      \\
                                      & =
                                      &\sqrt{\left(\frac{D}{H}\right)^{2}+F+\frac{4sT_{max}}{r_{T}}H}
\end{eqnarray*}
with
$$F=\frac{4AqT_{max}^{2}}{H^{2}r_{T}^{2}r_{I}^{2}}
\left(A(\delta-r_{I})-d_{I}(r_{I}-r_{T})-r_{I}(q-r_{I}
-r_{T})+r_{T}q\right).$$ It follows that :
$$T^{\ast}=\frac{1}{2}\left(-\frac{D}{H}+\sqrt{\left
(\frac{D}{H}\right)^{2}+F+\frac{4sT_{max}}{r_{T}}H}\right).$$
\end{preuve}\\
The combination of the proposition~\ref{p3}, proposition~\ref{p1}
and the lemma~\ref{l1} leads to the following theorem :
\begin{theorem}
The model (\ref{s}) admits a unique infected equilibrium
$E^{\ast}=(T^{\ast},I^{\ast},V^{\ast})$ if and only if
$\mathcal{R}_{0}>1$, where
\begin{eqnarray*}
 T^{\ast}  &=&\frac{1}{2}\left(-\frac{D}{H}+\sqrt{\Big{(}
 \frac{D}{H}\Big{)}^{2}+F+\frac{4sT_{max}}{r_{T}}H}\right),  \\
 I^{\ast}  &=& T^{\ast}\left(\frac{A}{r_{I}}-1\right)+T_{max}\left(1-\frac{\delta}{r_{I}}\right), \\
 V^{\ast}  &=& \frac{(1-\varepsilon)pI^{\ast}}{c};
\end{eqnarray*}
where where :
$$D=AT_{max}\left(\frac{1}{r_{T}}\Big{(}1+\frac{d_{T}+q}{A}\Big{)}-\frac{\delta}{r_{I}}\Big{(}\frac{1}{r_{T}}+\frac{1}{A}\Big{)}-\frac{q}{r_{T}r_{I}}\right);$$
$$F=\frac{4AqT_{max}^{2}}{H^{2}r_{T}^{2}r_{I}^{2}}\Big{(}A(\delta-r_{I})-d_{I}(r_{I}-r_{T})-r_{I}(q-r_{I}-r_{T})+r_{T}q\Big{)};$$
$$H=\frac{A^{2}}{r_{I}r_{T}}+\frac{A}{r_{I}}-\frac{A}{r_{T}}$$
 and
$$A=\frac{(1-\theta)\beta pT_{max}}{c}.$$

When $\mathcal{R}_{0}\leq1$ the unique equilibrium is the uninfected
equilibrium point or the infection-free steady state
$E^{0}=(T^{0},0,0)$.
\end{theorem}
\subsection{Local stability analysis of the
model \ref{s1} at the equilibrium points} For the study of local
stability of the model (\ref{s}) at the equilibrium points, let us
consider once more the functions $f_{1}$, $f_{2}$ et $f_{3}$ given
by (\ref{f1}), (\ref{f2}) and (\ref{f1}) respectively.
\subsubsection{Case of the uninfected  equilibrium point
or infection-free steady state}
\begin{theorem}\label{as}
The infection-free steady state $E^{0}=(T^{0},0,0)$ of model
(\ref{s}) is locally asymptotically stable if $\mathcal{R}_{0} \leq
1$ and unstable if $\mathcal{R}_{0} > 1$.
\end{theorem}
\begin{preuve} The Jacobian matrix $J(E^{0})$
of the system (\ref{s}) at $E^{0}$ is as the following :\\
\begin{eqnarray*}
 J(E^{0})&=&\left(
   \begin{array}{ccc}
 \frac{\partial f_{1}}{\partial T}(E^{0}) & \frac{\partial f_{1}}{\partial I}(E^{0})  & \frac{\partial f_{1}}{\partial V}(E^{0})\\\\
  \frac{\partial f_{2}}{\partial T}(E^{0}) & \frac{\partial f_{2}}{\partial I}(E^{0})  &  \frac{\partial f_{2}}{\partial V}(E^{0})\\ \\
 \frac{\partial f_{3}}{\partial T}(E^{0}) & \frac{\partial f_{3}}{\partial I}(E^{0})  & \frac{\partial f_{3}}{\partial
 V}(E^{0})
  \end{array}
  \right)\\
  &&
\end{eqnarray*}
i.e.
\begin{eqnarray*}
J(E^{0})&=&\left(
   \begin{array}{ccc}
 r_{T}\Big{(}1-\frac{2T^{0}}{T_{max}}\Big{)}-d_{T} & -\frac{r_{T}T^{0}}{T_{max}}+q & -(1-\eta)\beta T^{0}\\\\
  0 & r_{I}\Big{(}1-\frac{T^{0}}{T_{max}}\Big{)}-\delta  &  (1-\eta)\beta T^{0} \\ \\
  0 & (1-\varepsilon)p & -c
  \end{array}
  \right).\\
  &&
\end{eqnarray*}
 Since
  $$ r_{T}-d_{T}=-\frac{s}{T^{0}}+\frac{r_{T}T^{0}}{T_{max}},$$
  it follows that :
\begin{eqnarray*}
 J(E^{0})&=&\left(
   \begin{array}{ccc}
 -\frac{s}{T^{0}}-\frac{r_{T}T^{0}}{T_{max}} & -\frac{r_{T}T^{0}}{T_{max}}+q & -(1-\eta)\beta T^{0}\\\\
  0 & r_{I}\Big{(}1-\frac{T^{0}}{T_{max}}\Big{)}-\delta  &  (1-\eta)\beta T^{0} \\\\
  0 & (1-\varepsilon)p & -c
  \end{array}
  \right).
\end{eqnarray*}
Now let us show that the eigenvalues of the matrix $J(E^{0})$ have
negative real part if and only if  $\mathcal{R}_{0}<1$. \\
Considering the expression of  $J(E^{0})$,
$-\frac{s}{T^{0}}-\frac{r_{T}T^{0}}{T_{max}}$ is a negative
eigenvalue of the matrix $J(E^{0})$. Now let us consider the
sub-matrix  $J_{1}(E^{0})$ defined by :
\begin{equation}
J_{1}(E^{0})=\left(
   \begin{array}{ccc}
  r_{I}\left(1-\frac{T^{0}}{T_{max}}\right)-\delta  &  (1-\eta)\beta T^{0} \\\\
  (1-\varepsilon)p & -c
\end{array}
  \right).
\end{equation}
The trace of  $J_{1}(E^{0})$ is :
\begin{eqnarray*}
Tr(J_{1}(E^{0})) & = &-c-\delta+r_{I}\left(1-\frac{T^{0}}{T_{max}}\right)\\
                 & = & \delta\left( -\frac{c}{\delta}-1+\frac{r_{I}}{\delta}\left(1-\frac{T^{0}}{T_{max}}\right)\right) \\
                 & = & \delta\left( -\frac{c}{\delta}-1+R_{0}-\frac{(1-\theta)}{c\delta}\beta pT^{0}\right) \\
                 & = & -c -\delta(1-R_{0})-\frac{(1-\theta)}{c}\beta
                 pT^{0},
\end{eqnarray*}
and the determinant of  $J_{1}(E^{0}$ is:
\begin{eqnarray*}
|J_{1}(E^{0})| & = &
-cr_{I}\left(1-\frac{T^{0}}{T_{max}}\right)+c\delta-(1-\theta)\beta
                 PT^{0}\\
                 & = &c\delta\left[-\frac{r_{I}}{\delta}\left(1-\frac{T^{0}}{T_{max}}\right)+1-\frac{(1-\theta)\beta PT^{0}}{c\delta}
                 \right]\\
                 & = & c\delta(1-\mathcal{R}_{0}).
\end{eqnarray*}
The system (\ref{s}) est locally asymptotically stable at
  $E^{0}$ if and only if
\begin{equation*}
-c -\delta(1-\mathcal{R}_{0})-\frac{(1-\theta)}{c}\beta pT^{0}<0,
\end{equation*}
and
\begin{equation*}
   c\delta(1-\mathcal{R}_{0})>0
\end{equation*}
i.e.
\begin{equation*}
  \mathcal{R}_{0}<1.
\end{equation*}
Therefore the model (\ref{s}), is locally asymptotically stable at
 $E^{0}=(T^{0},0,0)$ when
$\mathcal{R}_{0}<1$ and unstable when $\mathcal{R}_{0}>1.$ This
completes the proof of theorem~\ref{as}
\end{preuve}
\subsubsection{Case of infected equilibrium point}
We start this subsection by two preliminary lemmas.
\begin{lemma}\label{lc}
The Jacobian matrix $J(E^{\ast})$, of lemma~\ref{lc} ,  of  system
(\ref{s}) at $E^{\ast}$ is given by :
\begin{equation*}
 J(E^{\ast})=\left(
   \begin{array}{ccc}
  -\frac{s}{T^{\ast}}-\frac{r_{T}T^{\ast}}{T_{max}}-q\frac{I^{\ast}}{T^{\ast}} & -\frac{r_{T}T^{\ast}}{T_{max}}+q & -(1-q)\beta T^{\ast}
  \\
   & &\\\\
  -\frac{r_{I}I^{\ast}}{T_{max}}+(1-\eta)\beta V^{\ast} & -\frac{r_{I}I^{\ast}}{T_{max}}-\frac{(1-\eta)\beta V^{\ast}T^{\ast}}{I^{\ast}} & (1-\eta)\beta T^{\ast}
  \\ \\
0 & (1-\varepsilon)p & -c \\ \\
  \end{array}
  \right)
\end{equation*}
\end{lemma}
\begin{preuve}
See the Appendice for the proof.
\end{preuve}
\begin{lemma}\label{le}
The characteristic equation of the Jacobian matrix $J(E^{\ast})$ of
the system (\ref{s}) at $E^{\ast}$   is given
 by the following cubic equation :
$$\lambda^{3}+A_{1}\lambda^{2}+A_{2}\lambda+A_{3}=0;$$
where :
\begin{eqnarray*}
 A_{1}  &=& c+\frac{s}{T^{\ast}}+\frac{r_{T}T^{\ast}+r_{I}I^{\ast}
 +AT^{\ast}}{T_{max}}+q\frac{I^{\ast}}{T^{\ast}},  \\
 A_{2}  &=& \frac{cs}{T^{\ast}}+\frac{cr_{T}T^{\ast}+sA+cr_{I}
 I^{\ast}}{T_{max}}+q\frac{I^{\ast}}{T^{\ast}}(r_{I}-\delta)
+\frac{sr_{I}I^{\ast}}{T^{\ast}T_{max}}+\frac{r_{T}AT^{\ast}
(T^{\ast}+I^{\ast})}{T_{max}^{2}}+\frac{cqI^{\ast}}{T^{\ast}}
+q\frac{I^{\ast}}{T_{max}}, \\
 A_{3}  &=& \frac{csr_{I}I^{\ast}}{T^{\ast}
 T_{max}}+\frac{cA^{2}I^{\ast}T^{\ast}}
 {T_{max}^{2}}-\frac{cAr_{I}I^{\ast}T^{\ast}}{T_{max}^{2}}
+\frac{cAr_{T}I^{\ast}T^{\ast}}{T_{max}^{2}}
+q\frac{cI^{\ast}}{T^{\ast}}(r_{I}-\delta).
\end{eqnarray*}
\end{lemma}
\begin{preuve}
See the appendice for the  proof of lemma~\ref{le}.
\end{preuve}

Now let :
 $$\Delta_{2}=\begin{vmatrix}
A_{1}& 1\\
A_{3}  &A_{2}
  \end{vmatrix}$$
  \\
By  Routh-Hurwitz criteria\cite{murr}, we have the following
results.
\begin{theorem}
For model (\ref{s}), when $\mathcal{R}_{0} > 1$ is valid, the unique
endemic equilibrium $E^{\ast}$ is locally asymptotically stable if
$\Delta_{2} > 0 $ and unstable if $\Delta_{2} < 0$.
\end{theorem}
Especially, we have :
\begin{cor}\label{locstat2}
The infected steady state during the therapy $E^{\ast}$ of the model
(\ref{s})
 is locally asymptotically stable if $\mathcal{R}_{0} > 1$ and unstable if
$\mathcal{R}_{0} > 1$.
\end{cor}
\begin{preuve}
Since  $A_{1}>0$,  $A_{2}>0$, it remains to show that
$\Delta_{2}>0$.
\begin{eqnarray*}
\Delta_{2} & = &A_{1}A_{2}-A_{3}\\
 &= &
\Big{(}c+\frac{s}{T^{\ast}}+\frac{r_{T}T^{\ast}}{T_{max}}
+\frac{r_{I}I^{\ast}}{T_{max}}+\frac{AT^{\ast}}{T_{max}}+
q\frac{I^{\ast}}{T^{\ast}}\Big{)}\Big{(}\frac{cs}{T^{\ast}}
+\frac{sr_{I}I^{\ast}}{T^{\ast}T_{max}}+\frac{sA}{T_{max}}
+\frac{cr_{I}I^{\ast}}{T_{max}}+ +\frac{cr_{T}T^{\ast}}{T_{max}}\\\\
&& +\frac{r_{T}A(T^{\ast})^{2}}{T_{max}^{2}}+
  \frac{r_{T}AT^{\ast}I^{\ast}}{T_{max}^{2}}+cq\frac{I^{\ast}}{T^{\ast}}+q\frac{I^{\ast}}{T^{\ast}}(r_{I}-\delta)
+q\frac{I^{\ast}}{T_{max}}\Big{)}
 -\frac{csr_{I}I^{\ast}}{T^{\ast}T_{max}}-\frac{cA^{2}T^{\ast}I^{\ast}}{T_{max}^{2}}+\frac{cr_{I}AT^{\ast}I^{\ast}}{T_{max}^{2}}
 \\\\
 &&-\frac{cr_{I}AT^{\ast}I^{\ast}}{T_{max}^{2}} -qc\frac{I^{\ast}}{T^{\ast}}(r_{I}-\delta)\\\\
 &=&\frac{c^{2}s}{T^{\ast}}+\frac{csr_{I}I^{\ast}}
 {T^{\ast}T_{max}}+\frac{csA}{T_{max}}+\frac{c^{2}
 r_{I}I^{\ast}}{T_{max}}
 +\frac{c^{2}r_{T}T^{\ast}}{T_{max}}+\frac{cr_{T}
 A(T^{\ast})^{2}}{T_{max}^{2}}+\frac{cr_{T}AT^{\ast}
 I^{\ast}}{T_{max}^{2}}
 +c^{2}q\frac{I^{\ast}}{T^{\ast}}+cq\frac{I^{\ast}}{T^{\ast}}
 (r_{I}-\delta)\\\\
 && +cq\frac{I^{\ast}}{T_{max}}+
 \frac{cs^{2}}{(T^{\ast})^{2}}+\frac{s^{2}r_{I}I^{\ast}}
 {(T^{\ast})^{2}T_{max}}+\frac{s^{2}A}{T^{\ast}T_{max}}
 +\frac{csr_{I}I^{\ast}}{T^{\ast}T_{max}}+\frac{csr_{T}}
{T_{max}}+\frac{sr_{T}AT^{\ast}}{T_{max}^{2}}
+\frac{sr_{T}AI^{\ast}}{T_{max}^{2}}
+\frac{csqI^{\ast}}{(T^{\ast})^{2}}
\\\\
&& +\frac{sqI^{\ast}}{(T^{\ast})^{2}}(r_{I}-\delta)+
\frac{sqI^{\ast}}{T^{\ast}T_{max}}+\frac{csr_{T}}{T_{max}}
+\frac{sr_{I}r_{T}I^{\ast}}{T_{max}^{2}}+
\frac{sAr_{T}T^{\ast}}{T_{max}^{2}}+
\frac{cr_{I}r_{T}I^{\ast}T^{\ast}}{T_{max}^{2}}
+\frac{cr_{T}^{2}(T^{\ast})^{2}}{T_{max}^{2}}
\end{eqnarray*}
\begin{eqnarray*}
&&+\frac{r_{T}^{2}A(T^{\ast})^{3}}{T_{max}^{3}}
+\frac{r_{T}^{2}A(T^{\ast})^{2}I^{\ast}}{T_{max}^{3}}+
\frac{r_{T}^{2}A(T^{\ast})^{2}I^{\ast}}{T_{max}^{3}}
+cq\frac{r_{T}I^{\ast}}{T_{max}}+q\frac{r_{T}I^{\ast}}{T_{max}}
(r_{I}-\delta)\frac{qr_{T}I^{\ast}T^{\ast}}{T_{max}^{2}}
+\frac{csr_{I}I^{\ast}}{T^{\ast}T_{max}}
\\\\
&&+\frac{sr_{I}^{2}(I^{\ast})^{2}}{T^{\ast}T_{max}}
+\frac{sAr_{I}I^{\ast}}{T_{max}^{2}}+
\frac{cr_{I}^{2}(I^{\ast})^{2}}{T_{max}^{2}}
+\frac{cr_{I}r_{T}I^{\ast}T^{\ast}}{T_{max}^{2}}
+\frac{r_{I}r_{T}A(T^{\ast})^{2}I^{\ast}}{T_{max}^{3}}
+\frac{r_{I}r_{T}AT^{\ast}(I^{\ast})^{2}}{T_{max}^{3}}
+qc\frac{r_{I}I^{\ast}}{T^{\ast}T_{max}}
\\\\
&&+\frac{qr_{I}(I^{\ast})^{2}}{T^{\ast}T_{max}}(r_{I}-\delta)
+\frac{qr_{I}(I^{\ast})^{2}}{T_{max}^{2}}+\frac{csA}{T_{max}}
+\frac{sr_{I}I^{\ast}}{T_{max}^{2}}+
\frac{sA^{2}T^{\ast}}{T_{max}^{2}}
+\frac{cr_{I}AI^{\ast}T^{\ast}}{T_{max}}
+\frac{cr_{T}A(T^{\ast})^{2}}{T_{max}^{2}}\\\\
&&+\frac{r_{T}A^{2}(T^{\ast})^{3}}{T_{max}^{3}}
+\frac{r_{T}A^{2}(T^{\ast})^{2}I}{T_{max}^{3}}
+q\frac{cAI^{\ast}}{T_{max}}+q\frac{AI^{\ast}}{T_{max}}(r_{I}
-\delta)+q\frac{AI^{\ast}T^{\ast}}{T_{max}^{2}}
+q\frac{csI^{\ast}}{(T^{\ast})^{2}}
\end{eqnarray*}
\begin{eqnarray*}
&&+q\frac{sr_{I}(I^{\ast})^{2}}{(T^{\ast})^{2}T_{max}}
+q\frac{sAI^{\ast}}{T^{\ast}T_{max}}
+q\frac{cr_{I}(I^{\ast})^{2}}{T^{\ast}T_{max}}
+q\frac{cr_{T}I^{\ast}}{T_{max}}+
q\frac{r_{T}AI^{\ast}T^{\ast}}{T_{max}^{2}}
+q\frac{r_{T}A(I^{\ast})^{2}}{T_{max}^{2}}\\\\
&&+q^{2}\frac{c(I^{\ast})^{2}}{(T^{\ast})^{2}}
+q^{2}\frac{(I^{\ast})^{2}}{(T^{\ast})^{2}}(r_{I}-\delta)
+q^{2}\frac{(I^{\ast})^{2}}{T^{\ast}T_{max}}
-\frac{csr_{I}I^{\ast}}{T^{\ast}T_{max}}
-\frac{cA^{2}T^{\ast}I^{\ast}}{T_{max}^{2}}\\\\
&&+\frac{cr_{I}AT^{\ast}I^{\ast}}{T_{max}^{2}}-\frac{cr_{T}AT^{\ast}I^{\ast}}{T_{max}^{2}}-qc\frac{I^{\ast}}{T^{\ast}}(r_{I}-\delta)\\
&=&B+D
\end{eqnarray*}
where,
\begin{eqnarray*}
B&=&\frac{c^{2}s}{T^{\ast}}+\frac{csr_{I}I^{\ast}}{T^{\ast}T_{max}}
+\frac{c^{2}r_{I}I^{\ast}}{T_{max}}
+\frac{c^{2}r_{T}T^{\ast}}{T_{max}}
+\frac{cr_{T}A(T^{\ast})^{2}}{T_{max}^{2}}
+\frac{cr_{T}AT^{\ast}I^{\ast}}{T_{max}^{2}}\\\\
&&+c^{2}q\frac{I^{\ast}}{T^{\ast}}+cq\frac{I^{\ast}}{T^{\ast}}
(r_{I}-\delta)+cq\frac{I^{\ast}}{T_{max}}+\frac{cs^{2}}
{(T^{\ast})^{2}}+
\frac{s^{2}r_{I}I^{\ast}}{(T^{\ast})^{2}T_{max}}
+\frac{s^{2}A}{T^{\ast}T_{max}}\\\\
&&+\frac{csr_{I}I^{\ast}}{T^{\ast}T_{max}}
+\frac{sr_{T}AT^{\ast}}{T_{max}^{2}}
+\frac{sr_{T}AI^{\ast}}{T_{max}^{2}}
+\frac{csqI^{\ast}}{(T^{\ast})^{2}}+
\frac{sqI^{\ast}}{(T^{\ast})^{2}}(r_{I}-\delta)
+\frac{sqI^{\ast}}{T^{\ast}T_{max}}\\\\
&&+\frac{sr_{I}r_{T}I^{\ast}}{T_{max}^{2}}
+\frac{sAr_{T}T^{\ast}}{T_{max}^{2}}
+\frac{cr_{I}r_{T}I^{\ast}T^{\ast}}{T_{max}^{2}}
+\frac{cr_{T}^{2}(T^{\ast})^{2}}{T_{max}^{2}}
+\frac{r_{T}^{2}A(T^{\ast})^{3}}{T_{max}^{3}}
+\frac{r_{T}^{2}A(T^{\ast})^{2}I^{\ast}}{T_{max}^{3}}\\\\
&&+\frac{r_{T}^{2}A(T^{\ast})^{2}I^{\ast}}{T_{max}^{3}}
+cq\frac{r_{T}I^{\ast}}{T_{max}}+q\frac{r_{T}I^{\ast}}{T_{max}}
(r_{I}-\delta)+
\frac{qr_{T}I^{\ast}T^{\ast}}{T_{max}^{2}}+\frac{csr_{I}I^{\ast}}
{T^{\ast}T_{max}}+\frac{sr_{I}^{2}(I^{\ast})^{2}}{T^{\ast}T_{max}}
\\\\
&&+\frac{sAr_{I}I^{\ast}}{T_{max}^{2}}
+\frac{cr_{I}^{2}(I^{\ast})^{2}}{T_{max}^{2}}+\frac{cr_{I}r_{T}I^{\ast}T^{\ast}}{T_{max}^{2}}+
\frac{r_{I}r_{T}A(T^{\ast})^{2}I^{\ast}}{T_{max}^{3}}
+\frac{r_{I}r_{T}AT^{\ast}(I^{\ast})^{2}}{T_{max}^{3}}\\\\
&&+qc\frac{r_{I}I^{\ast}}{T^{\ast}T_{max}}+
\frac{qr_{I}(I^{\ast})^{2}}{T^{\ast}T_{max}}(r_{I}-\delta)+\frac{qr_{I}(I^{\ast})^{2}}{T_{max}^{2}}+
\frac{sr_{I}I^{\ast}}{T_{max}^{2}}+
\frac{sA^{2}T^{\ast}}{T_{max}^{2}}+
\frac{cr_{I}AI^{\ast}T^{\ast}}{T_{max}}\\\\
&&+\frac{cr_{T}A(T^{\ast})^{2}}{T_{max}^{2}}
+\frac{r_{T}A^{2}(T^{\ast})^{3}}{T_{max}^{3}}
+\frac{r_{T}A^{2}(T^{\ast})^{2}I}{T_{max}^{3}}+
q\frac{AI^{\ast}}{T_{max}}(r_{I}-\delta)
+q\frac{AI^{\ast}T^{\ast}}{T_{max}^{2}}\\\\
&&+q\frac{csI^{\ast}}{(T^{\ast})^{2}}
+q\frac{sr_{I}(I^{\ast})^{2}}{(T^{\ast})^{2}T_{max}}
+q\frac{sAI^{\ast}}{T^{\ast}T_{max}}
+q\frac{cr_{I}(I^{\ast})^{2}}{T^{\ast}T_{max}}
+q\frac{cr_{T}I^{\ast}}{T_{max}}+q\frac{r_{T}AI^{\ast}T^{\ast}}
{T_{max}^{2}}\\\\
&&+q\frac{r_{T}A(I^{\ast})^{2}}{T_{max}^{2}}
+q^{2}\frac{c(I^{\ast})^{2}}{(T^{\ast})^{2}}
+q^{2}\frac{(I^{\ast})^{2}}{(T^{\ast})^{2}}(r_{I}-\delta)
+q^{2}\frac{(I^{\ast})^{2}}{T^{\ast}T_{max}}
-\frac{csr_{I}I^{\ast}}{T^{\ast}T_{max}}+
\frac{cr_{I}AT^{\ast}I^{\ast}}{T_{max}^{2}}-
\frac{cr_{T}AT^{\ast}I^{\ast}}{T_{max}^{2}}\\\\
&&-qc\frac{I^{\ast}}{T^{\ast}}(r_{I}-\delta),
\end{eqnarray*}
and
\begin{equation}\label{d}
D=2\frac{csA}{T_{max}}+2\frac{csr_{T}}{T_{max}}+q\frac{cAI^{\ast}}{T_{max}}-\frac{cA^{2}T^{\ast}I^{\ast}}{T_{max}^{2}}.
\end{equation}
Since $B>0$, it remains to show that $D>0$. \\
From (\ref{s1}) we have :\\
$$s+r_{T}T^{\ast}-\frac{r_{T}(T^{\ast})^{2}}{T_{max}}-\frac{r_{T}T^{\ast}I^{\ast}}{T_{max}}
-d_{T}T^{\ast}-(1-\eta)\beta V^{\ast}T^{\ast}+qI^{\ast}=0,$$
  yet
 $$V^{\ast}=\frac{(1-\varepsilon)pI^{\ast}}{c};$$ hence,
\begin{eqnarray*}
(1-\eta)\beta V^{\ast}T^{\ast} & = & \frac{(1-\theta)p\beta I^{\ast}T^{\ast}}{c}\\
                 & = & \frac{AI^{\ast}T^{\ast}}{T_{max}}.
\end{eqnarray*}
Thus,
$$\frac{AI^{\ast}T^{\ast}}{T_{max}}=s+r_{T}T^{\ast}-d_{T}T^{\ast}-\frac{r_{T}(T^{\ast})^{2}}{T_{max}}-\frac{r_{T}T^{\ast}I^{\ast}}{T_{max}}
+qI^{\ast}$$ it follows that ;
$$-\frac{cA^{2}I^{\ast}T^{\ast}}{T_{max}^{2}}=-\frac{csA}{T_{max}}-\frac{cAr_{T}T^{\ast}}{T_{max}}+\frac{cAd_{T}T^{\ast}}{T_{max}}
+\frac{cAr_{T}(T^{\ast})^{2}}{T_{max}^{2}}+\frac{cAr_{T}T^{\ast}I^{\ast}}{T_{max}^{2}}
-\frac{cAqI^{\ast}}{T_{max}}.$$ Reporting this previous expression
in (\ref{d}) yields :
\begin{eqnarray*}
D & =
&2\frac{csA}{T_{max}}+2\frac{csr_{T}}{T_{max}}+q\frac{cAI^{\ast}}{T_{max}}-\frac{csA}{T_{max}}-\frac{cAr_{T}T^{\ast}}{T_{max}}\\
&&+\frac{cAd_{T}T^{\ast}}{T_{max}}+\frac{cAr_{T}(T^{\ast})^{2}}{T_{max}^{2}}+\frac{cAr_{T}T^{\ast}I^{\ast}}{T_{max}^{2}}
-\frac{cAqI^{\ast}}{T_{max}}\\ \\
& =
&\frac{csA}{T_{max}}+2\frac{csr_{T}}{T_{max}}-\frac{cAr_{T}T^{\ast}}{T_{max}}+\frac{cAd_{T}T^{\ast}}{T_{max}}
+\frac{cAr_{T}(T^{\ast})^{2}}{T_{max}^{2}}+\frac{cAr_{T}T^{\ast}I^{\ast}}{T_{max}^{2}}
\end{eqnarray*}
Taking especially $s=d_{T}T_{max}$ and $\delta=d_{T}$, we obtain :
$T^{\ast}=\frac{d_{T}T_{max}}{A}$. Thus,\\
$$D=\frac{csA}{T_{max}}+c\delta
r_{T}+\frac{cAd_{T}T^{\ast}}{T_{max}}+\frac{cAr_{T}(T^{\ast})^{2}}
{T_{max}^{2}}+\frac{cAr_{T}T^{\ast}I^{\ast}}{T_{max}^{2}},$$ and
$D>0$, therefore the system (\ref{s}) is locally  asymptotically
stable at $E^{\ast}.$ This completes the proof of
Corollary~\ref{locstat2}.
\end{preuve}

\subsection{ Global stability analysis of the
 system at equilibrium points}
  The global stability analysis of a dynamical system is usually a
very complex problem. One of the most efficient methods
 to
solve this problem is Lyapunov's theory. To build the functions of
Lyapunov we will follow the method proposed by  A. Korobeinikov
\cite{bull, lyap, koro}.
\subsubsection{Case of infection-free steady state}
\begin{theorem}\label{stabglo1}
The infection-free steady state $E^{0}=(T^{0},0,0)$ of the model
    (\ref{s}) is globally
    asymptotically stable if the basic reproduction
     number $R_{0}<1-\frac{q}{\delta}$ and unstable if
     $R_{0}> 1-\frac{q}{\delta}$.
\end{theorem}
\begin{preuve}
Consider the Lyapunov function :
    $$L(T,I,V)=T-T^{0}-T^{0}ln\frac{T}{T^{0}}+I+\frac{(1-\eta)\beta
T^{0}}{c}V.$$ $L_{1}$ is defined, continuous and positive definite
for all $T>0$, $I>0$, $V>0$. Also, the global minimum $L_{1}=0$
occurs at the infection free equilibrium $E^{0}$. Further,  function
$L_{1}$, along the solutions of system (1), satisfies :
\begin{eqnarray*}
 \frac{dL}{dt}& = &\frac{\partial L}{\partial T}\frac{dT}{dt}+\frac{\partial L}{\partial I}\frac{dI}{dt}+\frac{\partial L}{\partial V}\frac{dV}{dt},\\
                                      & = & \left(1-\frac{T^{0}}{T}\right)\dot{T}+\dot{I}+\frac{(1-\eta)\beta T^{0}}{c}\dot{V}, \\
                                      & = &(T-T^{0})\frac{\dot{T}}{T}+\dot{I}+\frac{(1-\eta)\beta T^{0}}{c}\dot{V}, \\
                                      & = &(T-T^{0})\left(\frac{s}{T}+r_{T}-\frac{r_{T}(T+I)}{T_{max}}-dT-(1-\eta)\beta
                                      V+q\frac{I}{T}\right)+(1-\eta)\beta VT\\
                                      &   &+
                                      r_{I}I\left(1-\frac{T+I}{T_{max}}\right)-\delta
                                      I + \frac{(1-\theta)}{c}\beta
                                      pT^{0}T-(1-\eta) \beta
                                      T^{0}V,
                                      \end{eqnarray*}
                                      \begin{eqnarray*}
                                     & = &(T-T^{0})\left(\frac{s}{T}+r_{T}-dT-\frac{r_{T}(T+I)}{T_{max}}+q\frac{I}{T}\right)-T(1-\eta)\beta
                                      V+T^{0}(1-\eta)\beta V\\
                                      & &+(1-\eta)\beta
                                      VT +r_{I}I\left(1-\frac{T+I}{T_{max}}\right)-\delta
                                      I + \frac{(1-\theta)}{c}\beta
                                      pT^{0}I-(1-\eta) \beta
                                      T^{0}V,\\
                                      &=&(T-T^{0})\left(\frac{s}{T}+r_{T}-dT-\frac{r_{T}(T+I)}{T_{max}}\right)+qI-qI\frac{T^{0}}{T}+r_{I}I\left(1-\frac{T+I}{T_{max}}\right)
                                      +\frac{(1-\theta)}{c}\beta
                                      pT^{0}-\delta)I;
\end{eqnarray*}
yet $$r_{T}-d_{T}=\frac{rT^{0}}{T_{max}}-\frac{s}{T^{0}}$$ hence,
further collecting terms, we have :
\begin{eqnarray*}
 \frac{dL}{dt}& =
  &(T-T^{0})\left(\frac{s}{T}+\frac{rT^{0}}{T_{max}}-\frac{s}{T^{0}}-\frac{r_{T}(T+I)}{T_{max}}\right)+r_{I}I(1-\frac{T+I}{T_{max}})
                                      +\frac{(1-\theta)}{c}\beta pT^{0}-\delta)I,\\
                                      & = & (T-T^{0})\left(-\frac{s}{TT^{0}}(T-T^{0})-\frac{r_{T}}{T_{max}}(T-T^{0})-\frac{r_{T}I}{T_{max}}\right)+ r_{I}I\left(1-\frac{T+I}{T_{max}}\right)\\
                                      &  &+(\frac{(1-\theta)}{c}\beta pT^{0}-\delta)I+qI-qI\frac{T^{0}}{T},
\end{eqnarray*}
\begin{eqnarray*}
                                      & = &
                                      -\frac{s}{TT^{0}}(T-T^{0})^{2}-\frac{r_{T}}{T_{max}}\left((T-T^{0})^{2}+(T-T^{0})I\right)+r_{I}I-\frac{r_{I}IT}{T_{max}}-\frac{r_{I}I^{2}}{T_{max}}\\
                                      &  & +\left(\frac{(1-\theta)}{c}\beta
                                      pT^{0}-\delta\right)I+qI-qI\frac{T^{0}}{T},
\end{eqnarray*}
\begin{eqnarray*}
                                      &= &-\frac{s}{TT^{0}}(T-T^{0})^{2}-\frac{r_{T}}{T_{max}}\left((T-T^{0})^{2}+(T-T^{0})I+\frac{r_{I}IT}{r_{T}}+\frac{r_{I}I^{2}}{r_{T}}\right)+r_{I}I\\
                                      & &+\left(\frac{(1-\theta)}{c}\beta pT_{0}-\delta\right)I+qI-qI\frac{T^{0}}{T},\\
                                      &= &-\frac{s}{TT_{0}}(T-T^{0})^{2}-\frac{r_{T}}{T_{max}}\Bigl[(T-T^{0})^{2}+(T-T^{0})I+\frac{r_{I}IT}{r_{T}}+\frac{r_{I}I^{2}}{r_{T}}+\frac{r_{I}}{r_{T}}IT^{0}
                                      -\frac{r_{I}}{r_{T}}IT^{0}\Bigr]
\end{eqnarray*}
\begin{eqnarray*}
                                      &&+r_{I}I+\left(\frac{(1-\theta)}{c}\beta pT^{0}-\delta\right)I+qI-qI\frac{T^{0}}{T},\\
                                      &=&-\frac{s}{TT^{0}}(T-T^{0})^{2}-\frac{r_{T}}{T_{max}}\Bigl((T-T^{0})^{2}+(T-T^{0})I+\frac{r_{I}IT}{r_{T}}(T-T^{0})+\frac{r_{I}I^{2}}{r_{T}}
                                      +\frac{r_{I}}{r_{T}}IT^{0}\Bigr)\\
                                      &&+r_{I}I+\left(\frac{(1-\theta)}{c}\beta pT^{0}-\delta\right)I+qI-qI\frac{T^{0}}{T},
\end{eqnarray*}
                                      \begin{eqnarray*}
                                      &= &-\frac{s}{TT^{0}}(T-T^{0})^{2}-\frac{r_{T}}{T_{max}}(T+I-T^{0})(T+\frac{r_{I}}{r_{T}}I-T^{0})-\frac{r_{I}}{T_{max}}IT^{0}+r_{I}I\\
                                      & &+\left(\frac{(1-\theta)}{c}\beta pT^{0}-\delta\right)I+qI-qI\frac{T^{0}}{T},\\
                                      &= &-\frac{s}{TT^{0}}(T-T^{0})^{2}-\frac{r_{T}}{T_{max}}(T+I-T^{0})(T+\frac{r_{I}}{r_{T}}I-T^{0})+\delta I\Bigl(\frac{r_{I}}{\delta}
                                      -\frac{r_{I}T^{0}}{\delta T_{max}}
                                       \frac{(1-\theta)}{c\delta}\beta
                                       pT^{0}-1\Bigr)\\
                                        & &+qI-qI\frac{T^{0}}{T}.
\end{eqnarray*}
Furthermore , $$\mathcal{R}_{0}=\frac{(1-\theta)}{c\delta}\beta
pT^{0}+\frac{r_{I}}{\delta}\left(1-\frac{T^{0}}{T_{max}}\right),$$
hence
\begin{eqnarray*}
\frac{dL}{dt} & =
&-\frac{s}{TT^{0}}(T-T^{0})^{2}-\frac{r_{T}}{T_{max}}(T+I-T^{0})(T+\frac{r_{I}}{r_{T}}I-T^{0})-qI\frac{T^{0}}{T}
                                       +\delta I(\mathcal{R}_{0}-1)+qI,\\
                                      &=&-\frac{s}{TT^{0}}(T-T^{0})^{2}-\frac{r_{T}}{T_{max}}(T+I-T^{0})(T+\frac{r_{I}}{r_{T}}I-T^{0})
                                      -qI\frac{T^{0}}{T}
                                      + \delta
                                      I(\mathcal{R}_{0}-1+\frac{q}{\delta}).
\end{eqnarray*}
Since $r_{I}\leq r_{T}$ and
$\mathcal{R}_{0}<1-\frac{q}{\delta}$, we have $\frac{dL}{dt}\leq 0$\\
and $\frac{dL}{dt}=0$ if and only if $T=T^{0}$ and $I=0$
simultaneously.\\
Therefore, the largest compact invariant subset of the set
$$  M= \Big{\{} (T, I, V)\in \Omega : \frac{dL}{dt}= 0 \Big{\}} $$
is the singleton $\{E^{0}\}$. By the Lasalle invariance
principle\cite{khalil1}, the infection-free equilibrium is globally
asymptotically stable if $\mathcal{R}_{0}<1-\frac{q}{\delta}$. We
have seen previously that if $\mathcal{R}_{0} > 1 $, at least one of
the eigenvalues of the Jacobian matrix evaluated at $E^{0}$ has a
positive real part. Therefore, the infection-free equilibrium
$E^{0}$ is unstable when $\mathcal{R}_{0}>1$. This completes the
proof of the theorem.
\end{preuve}
\begin{remark}
The Lyapunov function defined in the proof of theorem~\ref{stabglo1}
has been obtained following the general giving by Korobonikov
\cite{lyap, koro, bull} for the dynamic virus fondamental model.
\end{remark}
\subsubsection{Case of infected equilibrium point}
We recall :
\begin{remark}\label{rm}
According to (\ref{se1}, \ref{se2}, \ref{se3}) the infected
equilibrium point $E^{\ast}$ verify :
\begin{eqnarray*}
r_{T}-d_{T}   &=& -\frac{s}{T^{\ast}}+(1-\eta)\beta
V^{\ast}+\frac{r_{T}}{T_{max}}(T^{\ast}+I^{\ast})
-q\frac{I^{\ast}}{T^{\ast}},
\end{eqnarray*}
\begin{eqnarray*}
r_{I}-\delta  &=& -\frac{(1-\eta)\beta
V^{\ast}T^{\ast}}{I^{\ast}}+\frac{r_{I}}{T_{max}}(T^{\ast}+I^{\ast}),
\end{eqnarray*}
\begin{eqnarray*}
 c  &=&\frac{(1-\varepsilon)pI^{\ast}}{V^{\ast}}.
\end{eqnarray*}
\end{remark}
\begin{theorem}
Suppose that $r_{I}=r_{T}$, $s=d_{T}T_{max}$ and $\delta=d_{T}$.\\
Then the infected steady state during therapy $E^{\ast}$ of model
(\ref{s}) is globally
    asymptotically stable as soon as it exists.
\end{theorem}
\begin{preuve}
Consider the Lyapunov function defined by :
\begin{equation*}
L(T,I,V)=T-T^{\ast}-T^{\ast}\ln\frac{T}{T^{\ast}}+I-I^{\ast}-I^{\ast}\ln\frac{I}{I^{\ast}}+\frac{(1-\eta)\beta
T^{\ast}V^{\ast}}{(1-\varepsilon)pI^{\ast}}\bigl(V-V^{\ast}
-V^{\ast}\ln\frac{V}{V^{\ast}}\bigr).
\end{equation*}
Let us show that $\frac{dL}{dt}\leq 0$ and
 $\frac{dL}{dt}= 0$ if and only if
  $T=T^{\ast}$, $I=I^{\ast}$, $V=V^{\ast}$ simultaneously.
  \\\indent
 The time
derivative of $L$ along the trajectories of system (\ref{s}) is :
\begin{eqnarray*}
\frac{dL}{dt}&=&\frac{\partial L}{\partial
T}\frac{dT}{dt}+\frac{\partial L}{\partial
I}\frac{dI}{dt}+\frac{\partial L}{\partial V}\frac{dV}{dt}
\end{eqnarray*}
\begin{eqnarray*}
&=&\left(1-\frac{T^{\ast}}{T}\right)\dot{T}+\left(1-\frac{I^{\ast}}{I}\right)\dot{I}+\frac{(1-\eta)\beta
T^{\ast}V^{\ast}}{(1-\varepsilon)pI^{\ast}}\left(1-\frac{V^{\ast}}{V}\right)\dot{V}
\end{eqnarray*}
\begin{eqnarray*}
&=&(T-T^{\ast})\frac{\dot{T}}{T}+(I-I^{\ast})\frac{\dot{I}}{I}+\frac{(1-\eta)\beta
T^{\ast}V^{\ast}}{(1-\varepsilon)pI^{\ast}}(V-V^{\ast})\frac{\dot{V}}{V}.
\end{eqnarray*}
Collecting terms, and canceling identical terms with opposite signs,
yields :
\begin{eqnarray}
\nonumber \frac{dL}{dt}   &=&
(T-T^{\ast})\left(\frac{s}{T}+r_{T}-\frac{r_{T}(T+I)}{T_{max}}
-d_{T}-(1-\eta)\beta V+q\frac{I}{T}\right)\\
\nonumber            &&+\frac{(1-\eta)\beta
T^{\ast}V^{\ast}}{(1-\varepsilon)pI^{\ast}}
\left(\frac{V-V^{\ast}}{V}\right)\Big{(}
 (1-\varepsilon)pI-cV\Big{)}
 \\
            &&+(I-I^{\ast})\left((1-\eta)\frac{\beta
VT}{I}+r_{I}\left(1-\frac{T+I}{T_{max}}\right)-\delta\right).
\label{e}
\end{eqnarray}
Reporting equalities of remark~\ref{rm} into (\ref{e}),
 we have :
\begin{eqnarray*}
\frac{dL}{dt}&=&(T-T^{\ast})\Bigl[\frac{s}{T}-\frac{s}{T^{\ast}}
+(1-\eta)\beta V^{\ast}+\frac{r_{T}}{T_{max}}(T^{\ast}+I^{\ast})
-q\frac{I^{\ast}}{T^{\ast}}-\frac{r_{T}(T+I)}{T_{max}}
-(1-\eta)\beta V+q\frac{I}{T}\Bigr]
\end{eqnarray*}
\begin{eqnarray*}
&& +(I-I^{\ast})\Bigl[(1-\eta)\frac{\beta
VT}{I}-r_{I}\frac{T}{T_{max}}-r_{I}\frac{I}{T_{max}}
-\frac{(1-\eta)\beta
V^{\ast}T^{\ast}}{I^{\ast}}+\frac{r_{I}}{T_{max}}(T^{\ast}+I^{\ast}
)\Bigr]+(1-\eta)\beta
\\&&\\
&&+\frac{T^{\ast}V^{\ast}}{(1-\varepsilon)pI^{\ast}}
\left(\frac{V-V^{\ast}}{V}\right)\left(
 (1-\varepsilon)pI-\frac{(1-\varepsilon)pI^{\ast}}
 {V^{\ast}}V\right)
 \end{eqnarray*}
\begin{eqnarray*}
 &=&-\frac{s}{TT^{\ast}}(T-T^{\ast})^{2}-\frac{r_{T}}{T_{max}}
 (T-T^{\ast})^{2}-\frac{r_{T}}{T_{max}}(T-T^{\ast})(I-I^{\ast})
 -(1-\eta)\beta(T-T^{\ast})(V-V^{\ast})\\
  &&\\
 &&+(1-\eta)\beta\Bigl[\Big{(}\frac{VT}{I}-\frac{V^{\ast}T^{\ast}}
 {I^{\ast}}\Big{)}(I-I^{\ast})
 +\frac{T^{\ast}V^{\ast}}{(1-\varepsilon)pI^{\ast}}
 \frac{(V-V^{\ast})}{V}\Bigl((1-\varepsilon)pI
 -\frac{(1-\varepsilon)}{V^{\ast}}pI^{\ast}V\Bigr)\Bigr]\\
 &&\\
 &&
-q\frac{I^{\ast}}{T^{\ast}}(T-T^{\ast})+q\frac{I}{T}(T-T^{\ast})-
 \frac{r_{I}}{T_{max}}(T-T^{\ast})(I-I^{\ast})
 -\frac{r_{I}}{T_{max}}(I-I^{\ast})^{2}
 \end{eqnarray*}
\begin{eqnarray*}
&=&-\frac{s}{TT^{\ast}}(T-T^{\ast})^{2}-\frac{r_{T}}{T_{max}}
(T-T^{\ast})^{2}-\frac{(r_{T}+r_{I})}{T_{max}}(T-T^{\ast})
(I-I^{\ast})-\frac{r_{I}}{T_{max}}(I-I^{\ast})^{2}
\\
&&\\&&
+(1-\eta)\beta\biggl[\Big{(}\frac{VT}{I}-\frac{V^{\ast}T^{\ast}}
{I^{\ast}}\Big{)}(I-I^{\ast})\\
&&\\&&-(T-T^{\ast})(V-V^{\ast})+\frac{T^{\ast}V^{\ast}}{(1-\varepsilon)pI^{\ast}}\frac{(V-V^{\ast})}{V}\left((1-\varepsilon)pI
 -\frac{(1-\varepsilon)}{V^{\ast}}pI^{\ast}V\right)\biggr]
 \\
&&\\
 &&-q\frac{1}{TT^{\ast}}\Bigl(T^{2}I^{\ast}+
 (T^{\ast})^{2}I-T^{\ast}TI-T^{\ast}TI^{\ast}\Bigr)
 \end{eqnarray*}
\begin{eqnarray*}
&=&-\frac{s}{TT^{\ast}}(T-T^{\ast})^{2}-\frac{1}{T_{max}}
(r_{T}T+r_{I}I-r_{T}T^{\ast}-r_{I}I^{\ast})(T+I-T^{\ast}-I^{\ast})
\end{eqnarray*}
\begin{eqnarray*}
&&+(1-\eta)\beta
T^{\ast}V^{\ast}\Bigl(\frac{VT}{V^{\ast}T^{\ast}}-
\frac{VTI^{\ast}}{IV^{\ast}T^{\ast}}-\frac{V^{\ast}T^{\ast}I}
{I^{\ast}V^{\ast}T^{\ast}}+
\frac{V^{\ast}T^{\ast}I^{\ast}}{V^{\ast}T^{\ast}I^{\ast}}
-\frac{TV}{T^{\ast}V^{\ast}}+\frac{TV^{\ast}}{T^{\ast}V^{\ast}}\\
&&\\
&&+\frac{T^{\ast}V}{T^{\ast}V^{\ast}}-\frac{T^{\ast}V^{\ast}}{T^{\ast}V^{\ast}}+\frac{IV}{I^{\ast}V}-\frac{I^{\ast}V^{2}}{I^{\ast}VV^{\ast}}
-\frac{V^{\ast}I}{I^{\ast}V}+\frac{V^{\ast}I^{\ast}V}{V^{\ast}
I^{\ast}V}\Bigr)\\
&&\\&&-q\frac{1}{TT^{\ast}}\Bigl((T-T^{\ast})^{2}I^{\ast}
+(T^{\ast})^{2}(I-I^{\ast})+TT^{\ast}(I^{\ast}-I)\Bigr)\\&&\\
&=&-\frac{s}{TT^{\ast}}(T-T^{\ast})^{2}-\frac{1}{T_{max}}(r_{T}T
+r_{I}I-r_{T}T^{\ast}-r_{I}I^{\ast})(T+I-T^{\ast}-I^{\ast})\\
&&\\&&+(1-\eta)\beta T^{\ast}V^{\ast}\Bigl(1+\frac{T}{T^{\ast}}
-\frac{VTI^{\ast}}{IV^{\ast}T^{\ast}}-
\frac{V^{\ast}I}{I^{\ast}V}\Bigr)-
q\frac{1}{TT^{\ast}}\big{(}(T-T^{\ast})^{2}I^{\ast}\\&&\\
&&+(T^{\ast})(I-I^{\ast})(T^{\ast}-T)\Bigr).
\end{eqnarray*}
Note that \\
$$1+\frac{T}{T^{\ast}}-\frac{VTI^{\ast}}{IV^{\ast}T^{\ast}}-\frac{V^{\ast}I}{I^{\ast}V}=\Big{(}3-\frac{T^{\ast}}{T}-\frac{VTI^{\ast}}{IV^{\ast}T^{\ast}}
-\frac{V^{\ast}I}{I^{\ast}V}\big{)}+\big{(}\frac{T}{T^{\ast}}+\frac{T^{\ast}}{T}-2\Big{)}$$
and
$$\big{(}\frac{T}{T^{\ast}}+\frac{T^{\ast}}{T}-2\big{)}
=\frac{(T-T^{\ast})^{2}}{TT^{\ast}}.$$ According to (\ref{se1}),
$$s=(1-\eta)\beta
T^{\ast}V^{\ast}+\Big{(}d_{T}-r_{T}+r_{T}\frac{(T^{\ast}+I^{\ast})}{T_{max}}\Big{)}T^{\ast}-qI^{\ast}.$$
furthermore,
$$T^{\ast}+I^{\ast}=T_{max};$$
hence,
$$s=(1-\eta)\beta T^{\ast}V^{\ast}+d_{T}T^{\ast}-qI^{\ast}.$$
By hypothesis, $r_{T}=r_{I}$ this leads to :
\begin{eqnarray*}
\frac{dL}{dt}&=&\Big{(}-\frac{d_{T}}{T}+q\frac{I^{\ast}}{TT^{\ast}
-\frac{(1-\eta)\beta
V^{\ast}}{T}}\Big{)}(T-T^{\ast})^{2}-\frac{r_{T}}{T_{max}}
(T+I-T^{\ast}-I^{\ast})^{2}\\ &&\\
&&+(1-\eta)\beta
V^{\ast}T^{\ast}\Big{(}3-\frac{T^{\ast}}{T}-\frac{VTI^{\ast}}
{IV^{\ast}T^{\ast}}
-\frac{V^{\ast}I}{I^{\ast}V}\Big{)}+(1-\eta)\beta
V^{\ast}T^{\ast}\frac{(T-T^{\ast})^{2}}{TT^{\ast}}\\
&&\\
&&-q\frac{1}{TT^{\ast}}\Big{(}(T-T^{\ast})^{2}I^{\ast}+T^{\ast}(I-I^{\ast})(T^{\ast}-T)\Big{)}
\end{eqnarray*}
\begin{eqnarray*}
&=&-\frac{d_{T}}{T}(T-T^{\ast})^{2}-\frac{r_{T}}{T_{max}}(T+I-T^{\ast}-I^{\ast})^{2}-q\frac{1}{T}(I-I^{\ast})(T^{\ast}-T)\\
 &&\\&&+(1-\eta)\beta
V^{\ast}T^{\ast}\big{(}3-\frac{T^{\ast}}{T}-\frac{VTI^{\ast}}{IV^{\ast}T^{\ast}}
-\frac{V^{\ast}I}{I^{\ast}V}\big{)}
\end{eqnarray*}
\begin{eqnarray*}
&=&-\frac{d_{T}}{T}(T-T^{\ast})^{2}-\frac{r_{T}}{T_{max}}
(T+I-T^{\ast}-I^{\ast})^{2}-q\frac{1}{T}(I-I^{\ast})(T^{\ast}-T)\\&&\\
&&+(1-\eta)\beta
V^{\ast}T^{\ast}\Big{(}3-\frac{(T^{\ast})^{2}II^{\ast}VV^{\ast}+(I^{\ast}VT)^{2}+(IV^{\ast})^{2}TT^{\ast}}{TT^{\ast}II^{\ast}VV^{\ast}}\Big{)}\\
&=&-\frac{d_{T}}{T}(T-T^{\ast})^{2}-\frac{r_{T}}{T_{max}}(T+I-T^{\ast}-I^{\ast})^{2}-q\frac{1}{T}(I-I^{\ast})(T^{\ast}-T)\\&&\\
&&+\frac{3(1-\eta)\beta
V^{\ast}T^{\ast}}{TT^{\ast}II^{\ast}VV^{\ast}}\left(TT^{\ast}II^{\ast}VV^{\ast}-\frac{1}{3}\Big{(}(T^{\ast})^{2}II^{\ast}VV^{\ast}+
(I^{\ast}VT)^{2}+(IV^{\ast})^{2}TT^{\ast}\Big{)}\right).\\
\end{eqnarray*}

Yet $$\frac{1}{3}\Big{(}(T^{\ast})^{2}II^{\ast}VV^{\ast}+
(I^{\ast}VT)^{2}+(IV^{\ast})^{2}TT^{\ast}\Big{)}\geq
TT^{\ast}II^{\ast}VV^{\ast}$$ since the geometric mean is less than
or equal to the arithmetic mean.\\
It should be noted that  $\frac{dL}{dt}\leq 0$ and $\frac{d
L_{2}}{dt} = 0 $ holds if and only if $(T, X, V)$ take the steady
states values $(T^{*}, X^{*}, V^{*})$ . Therefore the infected
equilibrium point $E^{*}$ is globally asymptotically stable. This
completes the proof of this theorem.
\end{preuve}
\subsection{Some numerical simulations}
Some numerical simulations have been done in the case $ R_{0}< 1 $
to confirm theoretical result obtain on global stability for the
uninfected equilibrium.\\\indent The following curves, obtained
using the Maple software, show the real-time evolution of uninfected
hepatocytes, infected hepatocytes and viral load. The values of the
parameters are taken in the parameter range defined by the
table~(\ref{tb}) and the initial conditions are $T_{0}=10^{3}$,
$I_{0}=2$ and $V_{0}=1$.
\subsubsection{Evolution in time of uninfected cells,
infected cells and viral load when $\mathcal{R}_{0}<1$}.
 Parameters values : $s=10$, $r_{T}=0.05$,
$r_{I}=0.112$, $d_{T}=0.001$, $d_{I}=0.1$, $T_{max}=10^{7}$,
$\beta=10^{-7}$, $\eta=10^{-7}$, $\varepsilon=10^{-8}$, $p=1$,
$q=0.5$, $c=2$. These parameters values yields : $T^{0}=4160020$ and
$\mathcal{R}_{0}=0.4556<1$.

\begin{figure}[!h]
\centering
\begin{subfigure}[b]{0.3\textwidth}
\includegraphics[angle=0,height=5cm,width=\textwidth]
{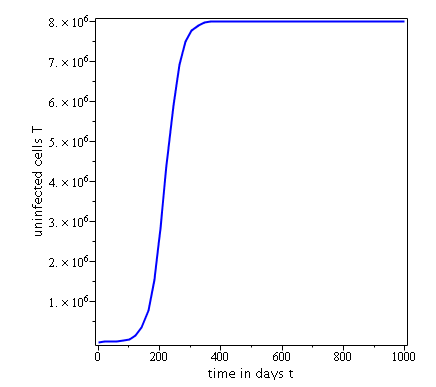} \subcaption{numerical solution curve for the uninfected
hepatocytes.}
\end{subfigure}
\begin{subfigure}[b]{0.3\textwidth}
\includegraphics[angle=0,height=5cm,width=\textwidth]{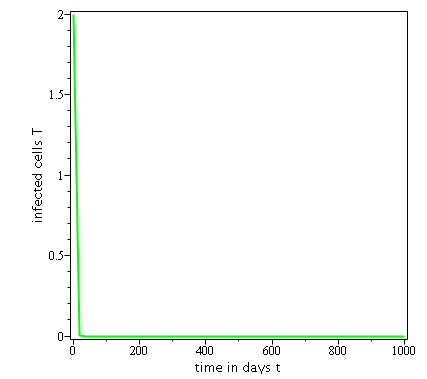}
\subcaption{numerical solution curve for the infected hepatocytes.}
\end{subfigure}
\begin{subfigure}[b]{0.3\textwidth}
\includegraphics[angle=0,height=5cm,width=\textwidth]{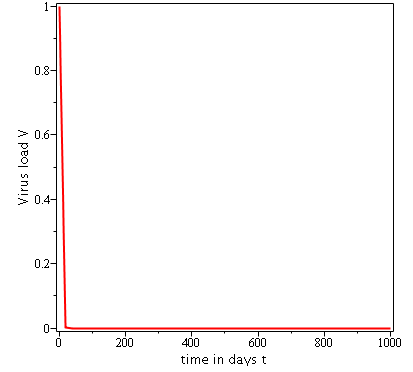}
\subcaption{numerical solution curve for the virus load.}
\end{subfigure}
\end{figure}
\begin{figure}[!h]
\centering
\includegraphics[angle=0,height=5cm,width=10cm]{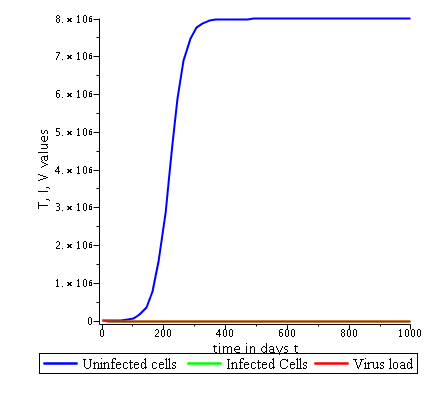}
\caption{Numerical simulations of the
 extended HCV model in 1000 days.}
\end{figure}

\newpage

 \subsubsection{Evolution in time of uninfected cells,
infected cells and viral load when $\mathcal{R}_{0}>1$.} Parameters
values: $s=10$, $r_{T}=2$, $r_{I}=0.112$, $d_{T}=0.01$, $d_{I}=0.3$,
$T_{max}=10^{7}$, $\beta=10^{-7}$, $\eta=10^{-4}$,
$\varepsilon=10^{-4}$, $p=1$, $q=0.5$, $c=0.5$. These parameters
values yields: $T^{0}=14875270$ et $\mathcal{R}_{0}=3.6501>1$.
\begin{figure}[!h]
\centering
\begin{subfigure}[b]{0.3\textwidth}
\includegraphics[angle=0,height=5cm,width=\textwidth]
{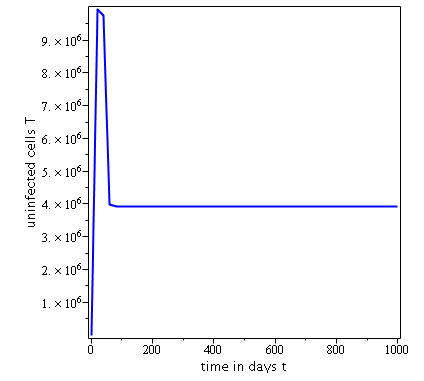} \subcaption{numerical solution curve for the uninfected
hepatocytes.}
\end{subfigure}
\begin{subfigure}[b]{0.3\textwidth}
\includegraphics[angle=0,height=5cm,width=\textwidth]{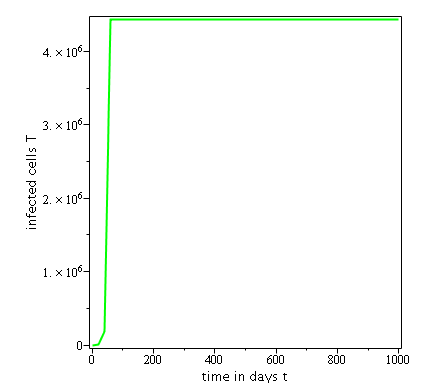}
\subcaption{numerical solution curve for the infected hepatocytes.}
\end{subfigure}
\begin{subfigure}[b]{0.3\textwidth}
\includegraphics[angle=0,height=5cm,width=\textwidth]{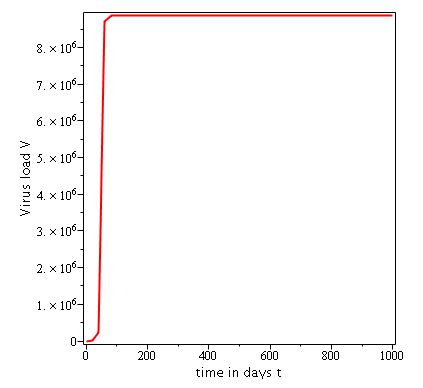}
\subcaption{numerical solution curve for the virus load.}
\end{subfigure}
\end{figure}
\begin{figure}[!h]
\centering
\includegraphics[angle=0,height=5cm,width=10cm]{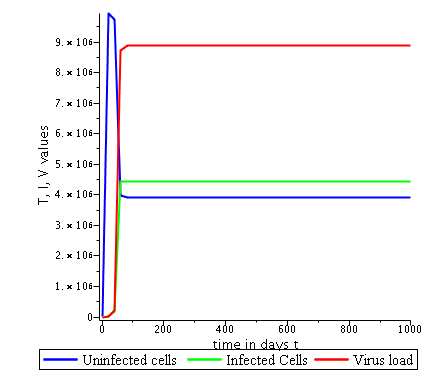}
\caption{Numerical simulations of the
 extended HCV model in 1000 days.}
\end{figure}

In short, in this section, it was a question to study  the global
stability of the model (\ref{s}). We have established that the model
(\ref{s}) is globally asymptotically stable at equilibrium points $
E^{\ast}$ and $ E^{0} $ when $ \mathcal{R}_{0} <1 $ and unstable
when $\mathcal{R}_{0}> 1 $. the numerical simulations has been
carried out using the Maple software confirming theoretical results.

\section*{Acknowledgement(s)}
I am  grateful to Professor Alan Rendall for valuable and tremendous
discussions. I wish to thank him for introducing me to Mathematical
Biology and to its relationship with Mathematical Analysis. I also
thank the Higher Teacher's Training College of the University of
Maroua were this paper were initiated.
\section*{Conclusion and discussions}
Having reached the end of our work, it emerges from all the
investigations presented that hepatitis C is a major health problem
in the world and especially in Cameroon. To understand the dynamics
of HCV and its infectious processes, mathematical models are present
as an important and unavoidable tool.   Global stability analysis
has been done, by the technique of Lyapunov, to the model of HCV
infection with proliferation cell and spontaneous healing, for
revealing  significant information for making  good decision for the
fighting against hepatitis C. We first show the existence of the
global solution to the Cauchy problem (\ref{s}), (\ref{s0}); then we
have calculated the basic reproduction ratio $\mathcal{R}_{0}$.
 We finally show that the only infected equilibrium point
 for the model is globally
asymptotically stable when $\mathcal{R}_{0} > 1 $ and unstable when
$ \mathcal{R}_{0} < 1$ under others hypotheses. Furthermore
uninfected equilibrium point
 for the model is globally
asymptotically stable when $\mathcal{R}_{0} < 1 $ and unstable
 when
$ \mathcal{R}_{0} > 1$. These theoretical results have been
confirmed by numerical simulations done using the software Maple.
Given the results obtained, this work is the beginning point of very
interesting other future investigations.\\
 We plan to extend our
analysis by focusing on more realistic models such as:
\begin{enumerate}
    \item models with  delay which involve delay ordinary differential equations;
    \item models taking into account space which involve Partial differential
equations;
    \item models taking into account  random phenomena which evolve stochastic
differential equations.
\end{enumerate}
 We also plan to focus on others methods of studying global
stability  like the  geometric method that can provides results with
less hypotheses on model (\ref{s}).

\section{Appendices}
\appendix
\section{Proof lemma~\ref{lc}}
\begin{preuve}
The Jacobian matrix $J(E^{\ast})$ of the system (\ref{s}) at
$E^{\ast}$ is given by :
\begin{eqnarray*}
 J(E^{\ast})&=&\left(
   \begin{array}{ccc}
 \frac{\partial f_{1}}{\partial T}(E^{\ast}) & \frac{\partial f_{1}}{\partial I}(E^{\ast})  & \frac{\partial f_{1}}{\partial V}(E^{\ast})\\\\
  \frac{\partial f_{2}}{\partial T}(E^{\ast}) & \frac{\partial f_{2}}{\partial I}(E^{\ast})  &  \frac{\partial f_{2}}{\partial V}(E^{\ast})\\ \\
 \frac{\partial f_{3}}{\partial T}(E^{\ast}) & \frac{\partial f_{3}}{\partial I}(E^{\ast})  & \frac{\partial f_{3}}{\partial
 V}(E^{\ast})
  \end{array}
  \right).\\
\end{eqnarray*}
Let us determine the coefficients of that Jacobian matrix .\\
$$\frac{\partial f_{1}}{\partial T}(E^{\ast})  =
r_{T}-\frac{r_{T}T^{\ast}}{T_{max}}-\frac{r_{T}I^{\ast}}{T_{max}}-\frac{r_{T}T^{\ast}}{T_{max}}-d_{T}-(1-\eta)\beta
V^{\ast}.$$ Yet
 $$V^{\ast}=\frac{(1-\varepsilon)}{c}pI^{\ast},$$
$$I^{\ast}=\left(\frac{(1-\theta)\beta
pT_{max}}{cr_{I}}-1\right)T^{\ast}+\frac{T_{max}}{r_{I}}(r_{I}-\delta)$$
and
\begin{eqnarray*}
   && -d_{T}=+\frac{(1-\theta)}{c}\beta p\left(\frac{r_{T}}{r_{I}}+\frac{(1-\theta)\beta
                 pT_{max}}{cr_{I}}-1\right)T^{\ast}
                 -\frac{s}{T^{\ast}}-\frac{qT_{max}}{T^{\ast}r_{I}}(r_{I}-\delta)  \\
   &&+\frac{r_{T}}{r_{I}}(r_{I}-\delta)
+\frac{(1-\theta)\beta
                 pT_{max}}{cr_{I}}(r_{I}-\delta)-\frac{(1-\theta)\beta pT_{max}}{cr_{I}}q-r_{T}+q;
\end{eqnarray*}
hence ,
\begin{eqnarray*}
\frac{\partial f_{1}}{\partial T}(E^{\ast})& = &
r_{T}-2\frac{r_{T}T^{\ast}}{T_{max}}-\frac{r_{T}}{T_{max}}\left(\Big{(}\frac{(1-\theta)\beta
                 pT_{max}}{cr_{I}}-1\Big{)}T^{\ast}+\frac{T_{max}}{r_{I}}(r_{I}-\delta)\right)\\
                 &&\\
                 &&+\frac{(1-\theta)}{c}\beta p\left(\frac{r_{T}}{r_{I}}+\frac{(1-\theta)\beta
                 pT_{max}}{cr_{I}}-1\right)T^{\ast}-\frac{s}{T^{\ast}}-\frac{qT_{max}}{T^{\ast}r_{I}}(r_{I}-\delta)
\end{eqnarray*}
\begin{eqnarray*}
                 &&+\frac{r_{T}}{r_{I}}(r_{I}-\delta)+\frac{(1-\theta)\beta
                 pT_{max}}{cr_{I}}(r_{I}-\delta)-\frac{(1-\theta)\beta pT_{max}}{cr_{I}}q-r_{T}\\
                 &&\\
                 &&+q-\frac{(1-\eta)\beta p(1-\varepsilon)}{c}\left(\Big{(}\frac{(1-\theta)\beta pT_{max}}{cr_{I}}-1\Big{)}T^{\ast}
                 +\frac{T_{max}}{r_{I}}(r_{I}-\delta)\right),\\
                 &&\\
                 &=&
                 -\frac{s}{T^{\ast}}+r_{T}-\frac{2r_{T}T^{\ast}}{T_{max}}+\frac{r_{T}T^{\ast}}{T_{max}}
                 -\frac{r_{T}}{cr_{I}}(1-\theta)\beta
                 pT^{\ast}-r_{T}+\frac{r_{T}}{r_{I}}\delta\\
                 &&\\
                 &&-\frac{(1-\theta)}{c}\beta pT^{\ast}+\frac{r_{T}}{cr_{I}}(1-\theta)\beta
                 pT^{\ast}+\frac{(1-\theta)^{2}\beta^{2}p^{2}}{c^{2}}\frac{T_{max}}{r_{I}}T^{\ast}-\frac{qT_{max}}{T^{\ast}}\\
                 &&\\
                 &&+\frac{qT_{max}}{T^{\ast}r_{I}}\delta+r_{T}-\frac{r_{T}}{r_{I}}\delta+\frac{(1-\theta)}{cr_{I}}\beta
                 pT_{max}(r_{I}-\delta)-\frac{(1-\theta)}{cr_{I}}\beta
                 pT_{max}q\\
                 &&\\
                 &&-r_{T}+q+\frac{(1-\theta)}{c}\beta
                 pT^{\ast}+\frac{(1-\theta)^{2}\beta^{2}p^{2}}{c^{2}}\frac{T_{max}}{r_{I}}T^{\ast}-\frac{(1-\theta)}{c}\beta
                 pT_{max}\\
                 &&\\
                 &&+\frac{(1-\theta)}{cr_{I}}\beta
                 pT_{max}\delta,\\
                 &&\\
                 &=&-\frac{s}{T^{\ast}}-\frac{r_{T}T^{\ast}}{T_{max}}-\frac{qI^{\ast}}{T^{\ast}};
\end{eqnarray*}
\begin{eqnarray*}
\frac{\partial f_{1}}{\partial I}(E^{\ast})   &=&
-\frac{r_{T}T^{\ast}}{T_{max}}+q;
\end{eqnarray*}

\begin{eqnarray*}
\frac{\partial f_{1}}{\partial V}(E^{\ast})   &=&-(1-\eta)\beta
T^{\ast} ;
\end{eqnarray*}

\begin{eqnarray*}
 \frac{\partial f_{2}}{\partial
T}(E^{\ast})  &=& -\frac{r_{I}I^{\ast}}{T_{max}}+(1-\eta)\beta
V^{\ast};
\end{eqnarray*}
\begin{eqnarray*}
\frac{\partial f_{2}}{\partial I}(E^{\ast}) & = &
r_{I}\big{(}1-\frac{T^{\ast}+I^{\ast}}{T_{max}}\big{)}
-\frac{r_{I}I^{\ast}}{T_{max}}-\delta,
\end{eqnarray*}
\begin{eqnarray*}
                 & = &
                 r_{I}-\frac{r_{I}T^{\ast}}{T_{max}}-2\frac{r_{I}I^{\ast}}{T_{max}}-\delta,
                \end{eqnarray*}
\begin{eqnarray*}
                 & = & -\frac{r_{I}I^{\ast}}{T_{max}}+r_{I}\left(1-\frac{T^{\ast}+\Big{(}\frac{(1-\theta)\beta
                 pT_{max}}{cr_{I}}-1\Big{)}T^{\ast}+\frac{T_{max}}{r_{I}(r_{I}-\delta)}}{T_{max}}\right)-\delta,
\end{eqnarray*}
\begin{eqnarray*}
                 &=&-\frac{r_{I}I^{\ast}}{T_{max}}+r_{I}-\frac{r_{I}T^{\ast}}{T_{max}}-\frac{(1-\theta)\beta
                 pT^{\ast}}{c}+\frac{r_{I}I^{\ast}}{T_{max}}-r_{I}+\delta-\delta,
\end{eqnarray*}
\begin{eqnarray*}
                 &=&-\frac{r_{I}I^{\ast}}{T_{max}}-\frac{(1-\theta)\beta
                 pT^{\ast}}{c},
                 \end{eqnarray*}
\begin{eqnarray*}
\frac{\partial f_{2}}{\partial I}(E^{\ast})
&=&-\frac{r_{I}I^{\ast}}{T_{max}}-\frac{(1-\eta)\beta
                 V^{\ast}T^{\ast}}{I^{\ast}}.
\end{eqnarray*}
\begin{eqnarray*}
 \frac{\partial f_{2}}{\partial V}(E^{\ast}) &=&(1-\eta)\beta
 T^{\ast};
\end{eqnarray*}
\begin{eqnarray*}
 \frac{\partial f_{3}}{\partial
T}(E^{\ast}) &=& 0 ;
\end{eqnarray*}
\begin{eqnarray*}
  \frac{\partial f_{3}}{\partial
I}(E^{\ast}) &=& (1-\varepsilon)p ;
\end{eqnarray*}

\begin{eqnarray*}
 \frac{\partial f_{3}}{\partial
V}(E^{\ast}) &=& -c .
\end{eqnarray*}
Therefore,
 $J(E^{\ast})=\left(
   \begin{array}{ccc}
  -\frac{s}{T^{\ast}}-\frac{r_{T}T^{\ast}}{T_{max}}
  -q\frac{I^{\ast}}{T^{\ast}} & -\frac{r_{T}T^{\ast}}{T_{max}}+q & -(1-q)\beta T^{\ast}
  \\
   & &\\\\
  -\frac{r_{I}I^{\ast}}{T_{max}}+(1-\eta)\beta V^{\ast} &
   -\frac{r_{I}I^{\ast}}{T_{max}}-\frac{(1-\eta)\beta V^{\ast}T^{\ast}}{I^{\ast}} & (1-\eta)\beta T^{\ast}
  \\ \\
0 & (1-\varepsilon)p & -C \\ \\
  \end{array}
  \right).$
  \\
  This completes the proof of the
lemma~\ref{lc}.
\end{preuve}

\section{Proof of lemma~\ref{le}}
\begin{preuve}
The characteristic equation is given by $|J-\lambda I|=0$,
i.e.\\
\begin{center}
  $\begin{vmatrix}
  -\frac{s}{T^{\ast}}-\frac{r_{T}T^{\ast}}{T_{max}}-q\frac{I^{\ast}}{T^{\ast}}-\lambda & -\frac{r_{T}T^{\ast}}{T_{max}}+q & -(1-q)\beta T^{\ast}
  \\
   & &\\\\
  -\frac{r_{I}I^{\ast}}{T_{max}}+(1-\eta)\beta V^{\ast} & -\frac{r_{I}I^{\ast}}{T_{max}}-\frac{(1-\eta)\beta V^{\ast}T^{\ast}}{I^{\ast}}-\lambda & (1-\eta)\beta T^{\ast}
  \\ \\
0 & (1-\varepsilon)p & -C-\lambda \\ \\
  \end{vmatrix}$=0.
\end{center}
With

$$q\frac{I^{\ast}}{T^{\ast}}=\frac{qT_{max}}{T^{\ast}}\left(-\frac{\delta}{r_{I}}+1\right)
  -q\left(\frac{A}{r_{I}}-1\right),$$
it follows that :
$$|J-\lambda I|=0$$
if and only if
$\left[-\frac{s}{T^{\ast}}-\frac{r_{T}T^{\ast}}{T_{max}}-\frac{qT_{max}}{T^{\ast}}\Big{(}-\frac{\delta}{r_{I}}+1\Big{)}
  -q\Big{(}\frac{A}{r_{I}}-1\Big{)}-\lambda\right]
  \Big{[}(c+\lambda)\Big{(}\frac{r_{I}I^{\ast}}{T_{max}}+\frac{(1-\eta)\beta
  V^{\ast}T^{\ast}}{I^{\ast}}
  +\lambda\Big{)}
  \\\\  -(1-\varepsilon)p(1-\eta)\beta
  T^{\ast}\Big{]}
  \\\\
   +\left[\frac{r_{I}I^{\ast}}{T_{max}}+(1-\eta)\beta
  V^{\ast}\right]
  \Bigl[(c+\lambda)\Big{(}\frac{r_{T}T^{\ast}}{T_{max}}-q\Big{)}+
  +(1-\varepsilon)p(1-\eta)\beta
  T^{\ast}\Bigr]=0,   $
  i.e.
\begin{equation*}
\Big{[}-\frac{s}{T^{\ast}}-\frac{r_{T}T^{\ast}}{T_{max}}-\frac{qT_{max}}{T^{\ast}}\Big{(}-\frac{\delta}{r_{I}}+1\Big{)}
  -q\Big{(}\frac{A}{r_{I}}-1\Big{)}-\lambda\Big{]}
  \Big{[}\lambda^{2}+\lambda(c+\frac{r_{I}I^{\ast}
  +AT^{\ast}}{T_{max}})+
  \frac{cr_{I}I^{\ast}}{T_{max}} \Big{]}
\end{equation*}
\begin{equation*}
 +\Big{[}\frac{r_{I}I^{\ast}-AI^{\ast}}{T_{max}}\Big{]}
  \Big{[}\lambda\Big{(}\frac{r_{T}T^{\ast}}{T_{max}}-q\Big{)}+\frac{cr_{T}T^{\ast}+AcT^{\ast}}{T_{max}}-qc\Big{]}=0,
\end{equation*}
  i.e
  \\\\
  $\lambda^{2}\Big{[}-\frac{s}{T^{\ast}}-\frac{r_{T}T^{\ast}}{T_{max}}-\frac{qT_{max}}{T^{\ast}}\Big{(}-\frac{\delta}{r_{I}}+1\Big{)}
  -q\Big{(}\frac{A}{r_{I}}-1\Big{)}\Big{]}+ \lambda\Big{[}
  \Big{(}c+\frac{r_{I}I^{\ast}+AT^{\ast}}{T_{max}}\Big{)}-\frac{s}{T^{\ast}}\\\\
  -\frac{r_{T}T^{\ast}}{T_{max}}-\frac{qT_{max}}{T^{\ast}}\Big{(}-\frac{\delta}{r_{I}}+1\Big{)}-q\Big{(}\frac{A}{r_{I}}
  -1)\Big{)}\Big{]}
  +\frac{cr_{I}I^{\ast}}{T_{max}}\Big{[}-\frac{s}{T^{\ast}}-\frac{r_{T}T^{\ast}}{T_{max}}-\frac{qT_{max}}{T^{\ast}}
  \Big{(}-\frac{\delta}{r_{I}}+1\Big{)}\\\\
  -q\Big{(}\frac{A}{r_{I}}-1\Big{)}\Big{]}
  -\lambda^{3}- \lambda^{2}\Big{(}c+\frac{r_{I}I^{\ast}+AT^{\ast}}{T_{max}}\Big{)}-\lambda\frac{cr_{T}I^{\ast}}{T_{max}}
  +\lambda\Big{(}\frac{r_{T}T^{\ast}}{T_{max}}-q\Big{)}\Big{(}\frac{r_{I}I^{\ast}-AI^{\ast}}{T_{max}}\Big{)}
  \\\\+\frac{r_{I}I^{\ast}-AI^{\ast}}{T_{max}}\Big{[}\frac{cr_{T}T^{\ast}+AcT^{\ast}}{T_{max}}-qc\Big{]}=0$.
\\\\
By developing the different factors in the previous equation, we get
: \\\\
$\lambda^{3}+\lambda^{2}\Big{[}\frac{s}{T^{\ast}}+\frac{r_{T}T^{\ast}}{T_{max}}+
\frac{qT_{max}}{T^{\ast}}\Big{(}-\frac{\delta}{r_{I}}+
1\Big{)}+q\Big{(}\frac{A}{r_{I}}-1\Big{)}+c+\frac{r_{I}I^{\ast}+AT^{\ast}}{T_{max}}\Big{]}\\\\
+\lambda\Big{[}(c+\frac{r_{I}I^{\ast}+AT^{\ast}}{T_{max}})(\frac{s}{T^{\ast}}+\frac{r_{T}T^{\ast}}{T_{max}}+
\frac{qT_{max}}{T^{\ast}}\Big{(}-\frac{\delta}{r_{I}}+1\Big{)}+q\Big{(}\frac{A}{r_{I}}-1)\Big{)}+\frac{cr_{I}I^{\ast}}{T_{max}}
+\Big{(}q-\frac{r_{T}T^{\ast}}{T_{max}}\Big{)}\Big{(}\frac{r_{I}I^{\ast}-AI^{\ast}}{T_{max}}\Big{)}\Big{]}
\\\\+
\Big{[}\frac{cr_{I}I^{\ast}}{T_{max}}\Big{(}\frac{s}{T^{\ast}}+\frac{r_{T}T^{\ast}}{T_{max}}+
\frac{qT_{max}}{T^{\ast}}\Big{(}-\frac{\delta}{r_{I}}+1\Big{)}q\Big{(}\frac{A}{r_{I}}-1\Big{)}\Big{)}+
\frac{r_{I}I^{\ast}-AI^{\ast}}{T_{max}}\Big{(}qc-\frac{cr_{T}T^{\ast}+AcT^{\ast}}{T_{max}}\Big{)}\Big{]}=0$.\\\\

Let:

$$A_{1}  =
\frac{s}{T^{\ast}}+\frac{r_{T}T^{\ast}}{T_{max}}+\frac{qT_{max}}{T^{\ast}}\Big{(}-\frac{\delta}{r_{I}}+1\Big{)}
                       +q\Big{(}\frac{A}{r_{I}}-1\Big{)}+c+
                       \frac{r_{I}I^{\ast}+AT^{\ast}}{T_{max}},$$
                       we have:
            $$A_{1}=c+\frac{s}{T^{\ast}}+\frac{r_{T}T^{\ast}+r_{I}I^{\ast}+AT^{\ast}}{T_{max}}+q\frac{I^{\ast}}{T^{\ast}}.$$
Let also :
\begin{eqnarray*}
 A_{2} &=& \Big{(}
c+\frac{r_{I}I^{\ast}+AT^{\ast}}{T_{max}}\Big{)}\left(\frac{s}{T^{\ast}}+\frac{r_{T}T^{\ast}}{T_{max}}+
\frac{qT_{max}}{T^{\ast}}(-\frac{\delta}{r_{I}}+1)+q(\frac{A}{r_{I}}-1)\right) \\
   && +\frac{cr_{I}I^{\ast}}{T_{max}}+
\Big{(}q-\frac{r_{T}T^{\ast}}{T_{max}}\Big{)}\Big{(}\frac{r_{I}I^{\ast}-AI^{\ast}}{T_{max}}\Big{)}
\end{eqnarray*}

we have :
\begin{eqnarray*}
A_{2}&=&\Big{(} c+\frac{r_{I}I^{\ast}}{T_{max}}+(1-\eta)\frac{\beta
V^{\ast}T^{\ast}}{I^{\ast}}\Big{)}\Big{(}\frac{s}{T^{\ast}}+\frac{r_{T}T^{\ast}}{T_{max}}+\frac{qT_{max}}
{T^{\ast}}(-\frac{\delta}{r_{I}}+1)
\end{eqnarray*}
\begin{eqnarray*}
 &&+q\Big{(}\frac{(1-\theta)}{cr_{I}}\beta
pT_{max}-1\Big{)}\Big{)}+\frac{cr_{I}I^{\ast}}{T_{max}}+\frac{c(1-\eta)\beta
V^{\ast}T^{\ast}}{I^{\ast}}-(1-\theta)\beta
pT^{\ast}\\
&&\\
&&+\Big{(}-\frac{r_{T}T^{\ast}}{T_{max}}+q\Big{)}\Big{(}\frac{r_{I}I^{\ast}}{T_{max}}-(1-\eta)\beta
V^{\ast}\Big{)}\\
&&\\
&=&\frac{cs}{T^{\ast}}+\frac{cr_{T}T^{\ast}}{T_{max}}+\frac{cqT_{max}}{r_{I}T^{\ast}}(r_{I}
-\delta)+cq\Big{(}\frac{(1-\theta)}{cr_{I}}\beta
pT_{max}-1\Big{)}+\frac{r_{I}I^{\ast}}{T_{max}}\frac{s}{T^{\ast}}\\
&&\\
&&+\frac{r_{I}I^{\ast}r_{T}T^{\ast}}{T_{max}^{2}}+q\frac{I^{\ast}}{T^{\ast}}(r_{I}-\delta)
+q\frac{r_{I}I^{\ast}}{T_{max}}\Big{(}\frac{(1-\theta)}{cr_{I}}\beta
pT_{max}-1\Big{)}\\
&&\\
&&+(1-\eta)\frac{\beta V^{\ast}s}{I^{\ast}}+\frac{(1-\eta)\beta
V^{\ast}(T^{\ast})^{2}r_{T}}{I^{\ast}T_{max}}+\frac{q(1-\eta)\beta
V^{\ast}T_{max}}{r_{I}I^{\ast}}(r_{I}-\delta)\\
&&\\
&&+\frac{q(1-\eta)\beta
V^{\ast}T^{\ast}}{I^{\ast}}\Big{(}\frac{(1-\theta)}{cr_{I}}\beta
pT_{max}-1\Big{)}+\frac{cr_{I}I^{\ast}}{T_{max}}+\frac{c(1-\eta)\beta
V^{\ast}T^{\ast}}{I^{\ast}}\\
&&\\
 &&-(1-\theta)\beta
pT^{\ast}-\frac{r_{T}T^{\ast}r_{I}I^{\ast}}{T_{max}^{2}}+(1-\eta)\beta
V^{\ast}\frac{r_{T}T^{\ast}}{T_{max}}+\frac{qr_{I}I^{\ast}}{T_{max}}-q(1-\eta)\beta
V^{\ast}
\end{eqnarray*}
\begin{eqnarray*}
&=&\frac{cs}{T^{\ast}}+\frac{cr_{T}T^{\ast}}{T_{max}}+\frac{cqT_{max}}{r_{I}T^{\ast}}(r_{I}-\delta)+cq\Big{(}\frac{(1-\theta)}{cr_{I}}\beta
pT_{max}-1\Big{)}+q\frac{I^{\ast}}{T^{\ast}}(r_{I}-\delta)\\
&&\\
&&+q\frac{r_{I}I^{\ast}}{T_{max}}\Big{(}\frac{(1-\theta)}{cr_{I}}\beta
pT_{max}-1\Big{)}+\frac{q(1-\eta)\beta
V^{\ast}T_{max}}{r_{I}I^{\ast}}(r_{I}-\delta)
\end{eqnarray*}
\begin{eqnarray*}
&&+\frac{q(1-\eta)\beta
V^{\ast}T^{\ast}}{I^{\ast}}\Big{(}\frac{(1-\theta)}{cr_{I}}\beta
pT_{max}-1\Big{)}+\frac{cr_{I}I^{\ast}}{T_{max}}+\frac{qr_{I}I^{\ast}}{T_{max}}-q(1-\eta)\beta
V^{\ast}\\
&&\\
&&+\frac{sr_{I}I^{\ast}}{T^{\ast}T_{max}}+\frac{sA}{T_{max}}+\frac{r_{T}AT^{\ast}(T^{\ast}+I^{\ast}}{T_{max}^{2}}\\
&&\\
&=&\frac{cs}{T^{\ast}}+\frac{cr_{T}T^{\ast}+sA+cr_{I}I^{\ast}}{T_{max}}+q\frac{I^{\ast}}{T^{\ast}}(r_{I}-\delta)
+\frac{sr_{I}I^{\ast}}{T^{\ast}T_{max}}+\frac{r_{T}AT^{\ast}(T^{\ast}+I^{\ast})}{T_{max}^{2}}+\frac{cqI^{\ast}}{T^{\ast}}+q\frac{I^{\ast}}{T_{max}}\\
\end{eqnarray*}

Let once more :  $$A_{3}
=\frac{cr_{I}I^{\ast}}{T_{max}}\Big{(}\frac{s}{T^{\ast}}+\frac{r_{T}T^{\ast}}{T_{max}}+
\frac{qT_{max}}{T^{\ast}}(-\frac{\delta}{r_{I}}+1)
+q(\frac{A}{r_{I}}-1)\Big{)}
+\frac{r_{I}I^{\ast}-AI^{\ast}}{T_{max}}
\Big{(}qc-\frac{cr_{T}T^{\ast} +AcT^{\ast}}{T_{max}}\Big{)}.$$ We
get :
\begin{eqnarray*}
A_{3}&=&\Big{(}\frac{cr_{I}I^{\ast}}{T_{max}}+\frac{c(1-\eta)\beta
V^{\ast}T^{\ast}}{I^{\ast}}-(-\theta)\beta
pT^{\ast}\Big{)}\Big{(}\frac{s}{T^{\ast}}+\frac{r_{T}T^{\ast}}{T_{max}}+
\frac{qT_{max}}{T^{\ast}}(-\frac{\delta}{r_{I}}+1)\\
&&\\
&&+q\Big{(}\frac{(1-\theta)\beta pT_{max}}{cr_{I}}-1\Big{)}
\Big{)}+\Big{(}-\frac{r_{I}I^{\ast}}{T_{max}}+(1-\eta)\beta
V^{\ast}\Big{)}\Big{(}\frac{cr_{T}T^{\ast}}{T_{max}}-cq+(1-\theta)\beta
pT^{\ast} \Big{)}\\
&&\\
&=&\frac{csr_{I}I^{\ast}}{T^{\ast}T_{max}}+\frac{cr_{I}r_{T}I^{\ast}T^{\ast}}{T_{max}^{2}}+
q\frac{cI^{\ast}}{T^{\ast}}(r_{I}-\delta)+q\frac{cr_{I}I^{\ast}}{T_{max}}\Big{(}\frac{(1-\theta)}{cr_{I}}\beta
pT_{max}-1\Big{)}\\
&&\\
&&+\frac{cs}{I^{\ast}}(1-\eta)\beta
V^{\ast}+\frac{cr_{T}(1-\eta)\beta
V^{\ast}(T^{\ast})^{2}}{I^{\ast}T_{max}}+q\frac{c(1-\eta)\beta
V^{\ast}}{r_{I}I^{\ast}}T_{max}(r_{I}-\delta)\\
&&\\
 &&+q\frac{c(1-\eta)\beta
V^{\ast}T^{\ast}}{I^{\ast}}\Big{(}\frac{(1-\theta)}{cr_{I}}\beta
pT_{max}-1\Big{)}-(1-\theta)s\beta
p-r_{T}(T^{\ast})^{2}\frac{(1-\theta)\beta p}{T_{max}}\\
&&\\
 &&-q\frac{(1-\theta)\beta
p}{r_{I}}T_{max}(r_{I}-\delta)-q(1-\theta)\beta
pT^{\ast}\Big{(}\frac{(1-\theta)}{cr_{I}}\beta
pT_{max}-1\Big{)}\\
&&\\
&&-\frac{cr_{I}r_{T}I^{\ast}T^{\ast}}{T_{max}^{2}}+cq\frac{r_{I}I^{\ast}}{T_{max}}-(1-\theta)\frac{r_{I}I^{\ast}T^{\ast}\beta
p}{T_{max}}+\frac{(1-\eta)\beta
V^{\ast}cr_{T}T^{\ast}}{T_{max}}\\
&&\\
&&-qc(1-\eta)\beta
V^{\ast}+(1-\theta)(1-\eta)\beta^{2}pV^{\ast}T^{\ast}\\
&&\\
&=&\frac{csr_{I}I^{\ast}}{T^{\ast}T_{max}}+\frac{cA^{2}I^{\ast}T^{\ast}}{T_{max}^{2}}-\frac{cAr_{I}I^{\ast}T^{\ast}}{T_{max}^{2}}
+\frac{cAr_{T}I^{\ast}T^{\ast}}{T_{max}^{2}}+q\frac{cI^{\ast}}{T^{\ast}}(r_{I}-\delta)\\
&&\\
&&+q\frac{cr_{I}I^{\ast}}{T_{max}}\Big{(}\frac{(1-\theta)}{cr_{I}}\beta
pT_{max}-1\Big{)}+q\frac{c(1-\eta)\beta
V^{\ast}}{r_{I}I^{\ast}}T_{max}(r_{I}-\delta)
\end{eqnarray*}
\begin{eqnarray*}
&&+q\frac{c(1-\eta)\beta
V^{\ast}T^{\ast}}{I^{\ast}}\Big{(}\frac{(1-\theta)}{cr_{I}}\beta
pT_{max}-1\Big{)}-q\frac{(1-\theta)\beta
p}{r_{I}}T_{max}(r_{I}-\delta)\\
&&\\
&&-q(1-\theta)\beta pT^{\ast}\Big{(}\frac{(1-\theta)}{cr_{I}}\beta
pT_{max}-1\Big{)}+cq\frac{r_{I}I^{\ast}}{T_{max}}-qc(1-\eta)\beta
V^{\ast}\\
&&\\
&=&\frac{csr_{I}I^{\ast}}{T^{\ast}T_{max}}+\frac{cA^{2}I^{\ast}T^{\ast}}{T_{max}^{2}}-\frac{cAr_{I}I^{\ast}T^{\ast}}{T_{max}^{2}}
+\frac{cAr_{T}I^{\ast}T^{\ast}}{T_{max}^{2}}+q\frac{cI^{\ast}}{T^{\ast}}(r_{I}-\delta).
\end{eqnarray*}
Therefore,
 $$\lambda^{3}+A_{1}\lambda^{2}+A_{2}\lambda+A_{3}=0;$$
 with
 $$A_{1}=c+\frac{s}{T^{\ast}}+\frac{r_{T}T^{\ast}+r_{I}I^{\ast}+AT^{\ast}}{T_{max}}+q\frac{I^{\ast}}{T^{\ast}}.$$

 $$A_{2}=\frac{cs}{T^{\ast}}+\frac{cr_{T}T^{\ast}+sA+cr_{I}I^{\ast}}{T_{max}}+q\frac{I^{\ast}}{T^{\ast}}(r_{I}-\delta)
+\frac{sr_{I}I^{\ast}}{T^{\ast}T_{max}}+\frac{r_{T}AT^{\ast}(T^{\ast}+I^{\ast})}{T_{max}^{2}}+\frac{cqI^{\ast}}{T^{\ast}}+q\frac{I^{\ast}}{T_{max}}$$
and
 $$A_{3}=\frac{csr_{I}I^{\ast}}{T^{\ast}T_{max}}+\frac{cA^{2}I^{\ast}T^{\ast}}{T_{max}^{2}}-\frac{cAr_{I}I^{\ast}T^{\ast}}{T_{max}^{2}}
+\frac{cAr_{T}I^{\ast}T^{\ast}}{T_{max}^{2}}+q\frac{cI^{\ast}}
{T^{\ast}}(r_{I}-\delta).$$ This completes the proof of the
lemma~\ref{le}.
\end{preuve}
\end{document}